\documentclass[10pt,leqno,twoside]{article}

\usepackage{amsmath,amssymb,amsthm}
\usepackage{titlesec,hyperref}
\usepackage{fancyhdr}

\usepackage{color}

\usepackage[margin=3.17cm]{geometry} %¿í°æʱÓã¬ÅäºÏË«±¶Ðоà

\usepackage{fancyhdr}

\titleformat{\subsection}{\it}{\thesubsection.\enspace}{1pt}{}

\newtheorem{theo}{Theorem}[section]
\newtheorem{lemm}[theo]{Lemma}
\newtheorem{defi}[theo]{Definition}
\newtheorem{coro}[theo]{Corollary}
\newtheorem{prop}[theo]{Proposition}
\newtheorem{rema}[theo]{Remark}
\numberwithin{equation}{section}

\def\ep{\varepsilon}

\allowdisplaybreaks
\begin{document}
\title{Well-posedness and global existence of 2D viscous shallow water system in Besov spaces
\hspace{-4mm}
}

\author{Yanan Liu$^1$
\quad Zhaoyang Yin$^2$ \\[10pt]
Department of Mathematics, Sun Yat-sen University,\\
510275, Guangzhou, P. R. China.\\[5pt]
}
\footnotetext[1]{Email: \it babyinarm@126.com}
\footnotetext[2]{Email: \it mcsyzy@mail.sysu.edu.cn}
\date{}
\maketitle

\begin{abstract}
In this paper we consider the Cauchy problem for 2D viscous
shallow water system in Besov spaces. We first establish the local
well-posedness of this problem in $B^s_{p,r}(\mathbb{R}^2)$,
$s>max\{1,\frac{2}{p}\}$ by using the Littlewood-Paley theory, the
Bony decomposition and the theories of transport equations and
transport diffusion equations. Then by the obtained local
well-posedness result, we can
prove the global existence of the system with small enough initial
data in $B^s_{p,r}(\mathbb{R}^2)$, $1\leq p\leq2$ and $s>\frac{2}{p}$.
Our obtained results improve the recent results in \cite{W}.

\vspace*{5pt}

\vspace*{5pt} \noindent{\it Keywords}: Viscous shallow water
system; Littlewood-Paley theory; Besov spaces; local
well-posedness; global existence
\end{abstract}

\vspace*{10pt}

%\phantomsection
%\addcontentsline{toc}{section}{\contentsname}
%Ìí¼ÓĿ¼µ½ÊéÇ©
\tableofcontents

\section{Introduction}
~~We consider the following Cauchy problems for 2D viscous shallow water equations
\begin{align}
\left\{
\begin{array}{l}
h(u_t+(u\cdot\nabla)u)-\nu\nabla\cdot(h\nabla u)+h\nabla h=0, \\[1ex]
h_t+div(hu)=0,  \\[1ex]
u|_{t=0}=u_0,\quad h|_{t=0}=h_0,
\end{array}
\right.
\end{align}
where $h(x,t)$ is the height of fluid surface,
$u(x,t)=(u^1(x,t),u^2(x,t))$ is the horizontal velocity field,
$x=(x_1,x_2)\in{\mathbb{R}}^2,$ and $0<\nu<1$ is the viscous
coefficients. For the initial data $h_0(x)$, we suppose that it is
a small perturbation of some positive constant $\bar h_0$. We
study the Cauchy problems (1.1) in Besov spaces
$B^s_{p,r}({\mathbb{R}}^2)$, $s>max\{\frac{2}{p},1\}$. For the
sake of convenient, we use the notation $B^s_{p,r}$ stands for
$B^s_{p,r}({\mathbb{R}}^2)$ in the following text. And we use the
notation $L^p$ stands for $L^p({\mathbb{R}}^2)$, and use the
notation $H^s$ stands for $H^s({\mathbb{R}}^2)$ as well.

Recently, Bresch et al. \cite{D, D.B} have systematically
introduced the viscous shallow water equations. Bui in \cite{bui}
proved the local existence and uniqueness of classical solutions
to the Cauchy-Dirichlet problem for the shallow water equations
with initial data in $C^{2+\alpha}$ by using Lagrangian
coordinates and H\"{o}lder space estimates. Kloeden in \cite{K}
and Sundbye in \cite{S1} independently showed the global existence
and uniqueness of classical solutions to the Cauchy-Dirichlet
problem using Sobolev space estimates by following the energy
method of Matsumura and Nishida \cite{M-N}. Sundbye in \cite{S2}
proved the existence and uniqueness of classical solutions to the
Cauchy problem using the method of \cite{M-N}. Wang and Xu in
\cite{W} obtained local solutions for any initial data and global
solutions for small initial data $h_0-\bar h_0,u_0\in H^s$, $s>2$.
Haspot got global existence in time for small initial data $h_0$,
$h_0-\bar h_0\in \dot{B}^{0}_{2,1}\cap \dot{B}^{1}_{2,1}$ and $u_0
\in \dot{B}^{0}_{2,1}$ as a special case in \cite{H}, and Chen,
Miao and Zhang in \cite{C-M-Z} to prove the local well-posedness with a more general diffusion
 in the space $u\in \mathcal{C}([0,T];B^0)\cap
L^1_T(B^2)$ and $h-\tilde{h}_0 \in \mathcal{C}([0,T];\tilde{B}^{0,1})\cap L^1(0,T;\tilde{B}^{2,1}),~h\geq \frac{1}{2}\tilde{h}_0$, with the initial data $(u_0,h_0-\bar{h}_0)\in B^0\times \tilde{B}^{0,1}, ~and~h_0\geq \bar{h}_0$. Moreover, they get the global existence with the small initial data.

There are two key points in this paper. In the proof of local well-posedness, we mainly use the Bony
decomposition to estimate the nonlinear terms. In the proof of the global existence with
small enough initial data in $B^s_{p,r}(\mathbb{R}^2)$, $1\leq p\leq2$
and $s>\frac{2}{p}$, we can't get the desired result after using once besov estimate for the transport-diffusion equation, the suitable iteration will be effective.

The main result of this paper is as follows:
\begin{theo}\label{t1}
Let $u_0,~h_0-\bar{h}_0\in B^s_{p,r}$, $s>max\{1,\frac{2}{p}\}$, $1\leq p\leq \infty$, $1\leq r<\infty$, $\|h_0-\bar{h}_0\|_{B^s_{p,r}}<<\bar{h}_0$. Then there exists a positive time $T$, a unique solution $(u,h)$ of the Cauchy problem (1.1) such that

$$u,h-\bar{h}_0\in \tilde{L}^{\infty}([0,T];B^s_{p,r})\cap C([0,T];B^s_{p,r}),~u\in \tilde{L}^2([0,T];B^{s+1}_{p,r}).$$
In the case $1\leq p\leq2,~~s>\frac{2}{p}$,$1\leq r<\infty$, there exists a small constant $\eta$ such that, if
$$\|h_0-\bar{h}_0\|_{B^s_{p,r}}+\|u_0\|_{B^s_{p,r}}\leq\eta,$$
the corresponding solution of the Cauchy problem (1.1)  is global in time.
\end{theo}

\section{Preliminaries}

First of all, we transform the system (1.1). For a sake of convenience, we take $\bar{h}_0=1$. Substituting $h$ by $1+h$ in (1.1), we have
\begin{align}
\left\{
\begin{array}{l}
u_t+(u\cdot\nabla)u-\nu\Delta u-\nu\nabla(\ln(1+h))\nabla u+\nabla h=0, \\[1ex]
h_t+divu+div(hu)=0,  \\[1ex]
u|_{t=0}=u_0,\quad h|_{t=0}=h_0,
\end{array}
\right.
\end{align}
here $h_0\in B^s_{p,r}$, and $\|h_0\|_{B^s_{p,r}}\leq \frac{1}{8C_0C_{s,p}}$, $C_0,C_{s,p}$ see Lemma \ref{l8}, \ref{l13} and \ref{l14}.
below.

Then we introduce the Littlewood-Paley decomposition briefly.

\begin{prop}\label{l2}
Littlewood-paley Decomposition:

\noindent Let $\mathcal{B}=\{\xi\in{\mathbb{R}^2},~
|\xi|\leq\frac{4}{3}\}$ be a ball, and
$\mathcal{C}=\{\xi\in{\mathbb{R}^2},
~\frac{3}{4}\leq|\xi|\leq\frac{8}{3} \}$ be an annulus. There
exist two radial functions $\chi$ and $\varphi$ valued in the
interval $[0,1]$, belonging respectively to
$\mathcal{D}(\mathcal{B})$ and $\mathcal{D}(\mathcal{C})$, such
that
\begin{align}
\forall~\xi\in{\mathbb{R}}^2,~~\chi(\xi)+\Sigma_{j\geq0}\varphi(2^{-j}\xi)=1,
\end{align}
\begin{align}
\forall~\xi\in{\mathbb{R}}^2\backslash\{0\},~~\Sigma_{j\in\mathbb{Z}}\varphi(2^{-j}\xi)=1,
\end{align}
\begin{align}
|j-j'|\geq2~\Rightarrow Supp\varphi(2^j\cdot)\cap Supp\varphi(2^{j'}\cdot)=\emptyset,
\end{align}
\begin{align}
j\geq2~\Rightarrow Supp\chi\cap Supp\varphi(2^j\cdot)=\emptyset,
\end{align}

the set $\tilde{\mathcal{C}}\overset{def}{=}~B(0,2/3)+\mathcal{C}$ is an annulus, and we have
\begin{align}
|j-j'|\geq5~\Rightarrow 2^j\tilde{\mathcal{C}}\cap 2^{j'}\mathcal{C}=\emptyset.
\end{align}

Further, we have
\begin{align}
\forall~\xi\in{\mathbb{R}}^2,~~\frac{1}{2}\leq\chi^2(\xi)+\Sigma_{j\geq0}\varphi^2(2^{-j}\xi)\leq1,
\end{align}
\begin{align}
\forall~\xi\in{\mathbb{R}}^2\backslash\{0\},~~\frac{1}{2}\leq\Sigma_{j\in\mathbb{Z}}\varphi^2(2^{-j}\xi)\leq1.
\end{align}
\end{prop}

Now we can define the nonhomogeneous dyadic blocks $\Delta_j$ and the nonhomogeneous low-frequency cut-off operator $S_j$ as follows:
$$\Delta_j u=0,~if~j\leq-2,~~\Delta_{-1}u=\chi(D)u=\int_{\mathbb{R}^2}\tilde{h}(y)u(x-y)dy,$$
$$\Delta_ju=\varphi(2^{-j}D)u=2^{jd}\int_{\mathbb{R}^2}h(2^jy)u(x-y)dy~~if~~j\geq0.$$
and
$$S_ju=\sum_{j'\leq j-1}\Delta_{j'}u.$$

Where $h=\mathcal{F}^{-1}\varphi$ and $\tilde{h}=\mathcal{F}^{-1}\chi$.

Next we define the Besov spaces:
\begin{defi}
Let $s\in R$ and $(p,r)\in[1,\infty]^2$. The nonhomogeneous Besov space $B^s_{p,r}$ consists of all tempered distribution $u$ such that:
$$\left(\sum_{j\geq-1}(2^{js}\|\Delta_ju\|_{L^p})\right)_{\ell^r}<\infty,$$
and naturally the Besov norm is defined as follows
$$\|u\|_{B^s_{p,r}}=\left(\sum_{j\geq-1}(2^{js}\|\Delta_ju\|_{L^p})\right)_{\ell^r}.$$
\end{defi}
\begin{defi}
The Bony decomposition:
The nonhomogeneous paraproduct of $v$ by $u$ is defined by
$$T_uv=\underset{j}{\sum}S_{j-1}u\Delta_jv.$$
The nonhomogeneous remainder of $u$ by $v$ is defined by
$$R(u,v)=\underset{|k-j|\leq1}{\sum}\Delta_ku\Delta_jv.$$
The operators $T$ and $R$ are bilinear, and we have the following Bony decomposition
$$uv=T_vu+T_uv+R(u,v).$$
\end{defi}

Then we give some properties of the Besov spaces which will be used in this paper.
\begin{lemm}\label{l5}
Let $1\leq p_1\leq p_2\leq\infty$ and $1\leq r_1\leq r_2\leq\infty$. Then for any $s\in \mathbb{R}$, the space
$B^s_{p_1,r_1}$ is continuously embedded in $B^{s-d(\frac{1}{p_1}-\frac{1}{p_2})}_{p_2,r_2}$. Obviously, we also have that the space $B^{s_2}_{p,r}$ is continuously embedded in $B^{s_1}_{p,r}$ and $B^{s_2}_{p,\infty}$ is continuously embedded in $B^{s_1}_{p,1}$ if $s_1<s_2$.
\end{lemm}
\begin{lemm}\label{l6}
If $u\in B^s_{p,r}$, then $\nabla u\in B^{s-1}_{p,r}$, and we have
$$\|\nabla u\|_{B^{s-1}_{p,r}}\leq C\|u\|_{B^s_{p,r}}.$$
\end{lemm}
\begin{lemm}\label{0}
If $s_1$ and $s_2$ are real numbers such that $s_1<s_2$,
$\theta\in(0,1)$, and $p,r\in[1,\infty]$, then we have
$$\|u\|_{B^{\theta s_1+(1-\theta)s_2}_{p,r}}\leq\|u\|^\theta_{B^{s_1}_{p,r}}\|u\|^{1-\theta}_{B^{s_2}_{p,r}}.$$
\end{lemm}
\begin{lemm}\label{l7}
The set $B^s_{p,r}$ is a Banach space and satisfies the Fatou property, namely, if $(u_n)_{n\in N}$ is a
bounded sequence of $B^s_{p,r}$. Then an element $u$ of $B^s_{p,r}$ and a subsequence $u_{\psi(n)}$ exist such that:

$\underset{n\rightarrow\infty}{lim}~u_{\psi(n)}=u$ in $\mathcal{S}'$,  $\|u\|_{B^s_{p,r}}\leq C \underset{n\rightarrow\infty}{lim} inf\|u_{\psi(n)}\|_{B^s_{p,r}}$.
\end{lemm}
\begin{lemm}\label{l8}
If $s>\frac{d}{p}$ or $s=\frac{d}{p},~r=1$, then the $B^s_{p,r}$ space is continuously embedded in $L^{\infty},$ i.e
$$\|u\|_{L^{\infty}}\leq C_{s,p}\|u\|_{B^s_{p,r}}.$$
\end{lemm}

\begin{lemm}\label{l9}
Let f be a smooth function, $f(0)=0$, $s>0,~(p,r)\in[1,\infty]^2$. If $u\in B^s_{p,r}\cap L^{\infty}$, then so does $f\circ u$, and we have
$$\|f\circ u\|_{B^s_{p,r}}\leq C\left(s,f',\|u\|_{L^{\infty}}\right)\|u\|_{B^s_{p,r}}.$$
\end{lemm}
\begin{lemm}\label{l10}
A constant $C$ exists which satisfies the following inequalities for any couple of real numbers $(s,t)$ with t negative and any $(p,r_1,r_2)$ in $[1,\infty]^3$:
$$\|T\|_{\mathcal{L}(L^{\infty}\times B^s_{p,r};B^s_{p,r})}\leq C^{|s|+1},$$
$$\|T\|_{\mathcal{L}(B^t_{\infty,r_1}\times B^s_{p,r_2};B^{s+t}_{p,r})}\leq \frac{C^{|s+t|+1}}{-t}
~~with~~\frac{1}{r}\overset{def}{=}min\{1,\frac{1}{r_1}+\frac{1}{r_2}\}.$$
\end{lemm}
\begin{lemm}\label{l11}
A constant $C$ exists which satisfies the following inequalities. Let $(s_1,s_2)$ be in $\mathbb{R}^2$ and
$(p_1,p_2,r_1,r_2)$ be in $[1,\infty]^4$. Assume that
$$\frac{1}{p}\overset{def}{=}\frac{1}{p_1}+\frac{1}{p_2}\leq1~~and~~
\frac{1}{r}\overset{def}{=}\frac{1}{r_1}+\frac{1}{r_2}\leq1.$$
If $s_1+s_2>0$, then we have, for any $(u,v)$ in $B^{s_1}_{p_1,r_1}\times B^{s_2}_{p_2,r_2}$,
$$\|R(u,v)\|_{B^{s_1+s_2}_{p,r}}\leq\frac{C^{|s_1+s_2|+1}}{s_1+s_2}\|u\|_{B^{s_1}_{p_1,r_1}}\|v\|_{B^{s_2}_{p_2,r_2}}.$$
If $r=1$ and $s_1+s_2=0$, then we have, for any $(u,v)$ in $B^{s_1}_{p_1,r_1}\times B^{s_2}_{p_2,r_2}$,
$$\|R(u,v)\|_{B^0_{p,\infty}}\leq C^{|s_1+s_2|+1}\|u\|_{B^{s_1}_{p_1,r_1}}\|v\|_{B^{s_2}_{p_2,r_2}}.$$
\end{lemm}
\begin{coro}\label{r}
Actually, the estimate of the remainder term $\|R(u,v)\|_{B^s_{p,r}}$ can be the same with $\|T_uv\|_{B^s_{p,r}}$ or $\|T_vu\|_{B^s_{p,r}}$. If $s>\frac{2}{p}-2$ in $p\in[1,2]$ or $s>-\frac{2}{p}$ in $p>2$.
\end{coro}
\begin{lemm}\label{l12}
For any $s>0$ and $(p,r)\in[1,\infty]^2$, the space $B^s_{p,r}\cap L^{\infty}$ is an algebra, and a constant exists such that:

$$\|uv\|_{B^s_{p,r}}\leq \frac{C^{s+1}}{s}\Big(\|u\|_{L^{\infty}}\|v\|_{B^s_{p,r}}+\|v\|_{L^{\infty}}\|u\|_{B^s_{p,r}}\Big).$$
Moreover, if $s>\frac{d}{p}~or~s=\frac{d}{p},r=1$, we have
$$\|uv\|_{B^s_{p,r}}\leq \frac{C^{s+1}}{s}\|u\|_{B^s_{p,r}}\|v\|_{B^s_{p,r}}.$$
\end{lemm}

For the transport equations
\begin{align}
\left\{
\begin{array}{l}
\partial_tf+v\cdot\nabla f=g\\
f_{|t=0}=f_0,
\end{array}
\right.
\end{align}
we have
\begin{lemm}\label{l13}
Let $1\leq p\leq p_1\leq\infty,~1\leq r\leq\infty$. Assume that
\begin{align}
s\geq-d\,min\left(\frac{1}{p_1},\frac{1}{p'}\right) \quad or \quad s\geq-1-d\,min\left(\frac{1}{p_1},\frac{1}{p'}\right)~if~div\,v=0
\end{align}
with strict inequality if $r<\infty$.

There exists a constant $C_0$, depending only on $d, p, p_1, r$ and $s$, such that for all solutions
$f\in L^{\infty}([0,T];B^s_{p,r})$ of (2.9), initial data $f_0$ in $B^s_{p,r}$, and $g$ in $L^1([0,T];B^s_{p,r})$, we have, for $a.e.\,t\in[0,T]$,
$$\|f\|_{\tilde{L}_t^{\infty}(B^s_{p,r})}\leq\left(\|f_0\|_{B^s_{p,r}}+
\int^t_0exp(-C_0V_{p_1}(t'))\|g(t')\|_{B^s_{p,r}}dt'\right)exp(C_0V_{p_1}(t))$$
with, if the inequality is strict in (2.10),
\begin{align}
V'_{p_1}(t)=\left\{\begin{array}{l}\|\nabla v(t)\|_{B^{s-1}_{p_1,r}},~if~s>1+\frac{d}{p_1}~or~s=1+\frac{d}{p_1},~r=1,\\
\|\nabla v(t)\|_{B^{\frac{d}{p_1}}_{p_1,\infty}\cap L^{\infty}},~if~s<1+\frac{d}{p_1}
\end{array}\right.
\end{align}
and, if equality holds in (2.10) and $r=\infty$,
$$V'_{p_1}=\|\nabla v(t)\|_{B^{\frac{d}{p_1}}_{p_1,1}}.$$
If $f=v$, then for all $s>0$($s>-1,~if~div\,u=0$), the estimate holds with
$$V'_{p_1}(t)=\|\nabla u\|_{L^{\infty}}$$.
Where
$\|u\|_{\tilde{L}^{\rho}_T(B^s_{p,r})}$ is defined in Lemma \ref{l15}.
\end{lemm}

For the transport diffusion equations
\begin{align}
\left\{
\begin{array}{l}
\partial_tf+v\cdot\nabla f-\nu\Delta f=g\\
f_{|t=0}=f_0,
\end{array}
\right.
\end{align}
we have the following lemma.
\begin{lemm}\label{l14}
Let $1\leq p\leq p_1\leq\infty,~1\leq r\leq\infty,~s\in\mathbb{R}$ satisfy (2.10), and $V_{p_1}$ be defined as in Lemma \ref{l13}.

There exists a constant $C_0$ which depends only on $d, r, s$ and $s-1-\frac{d}{p_1}$ and is such that for any smooth solution of (11) and $1\leq\rho_1\leq\rho\leq\infty,$ we have
$$\nu^{\frac{1}{\rho}}\|f\|_{\tilde{L}^{\rho}_T(B^{s+\frac{2}{\rho}}_{p,r})}\leq C_0e^{C_0(1+\nu T)^{\frac{1}{\rho}}V_{p_1}(T)}\Big((1+\nu T)^{\frac{1}{\rho}}\|f_0\|_{B^s_{p,r}}$$
$$+(1+\nu T)^{1+\frac{1}{\rho}-\frac{1}{\rho_1}}\nu^{\frac{1}{\rho_1}-1}\|g\|_
{\tilde{L}^{\rho_1}_T(B^{s-2+\frac{2}{\rho_1}}_{p,r})}\Big).$$
\end{lemm}

For the space $\tilde{L}^{\rho}_T(B^s_{p,r})$, we have the following properties:
\begin{lemm}\label{l15}
For all $T>0,~s\in\mathbb{R},$ and $1\leq r,\rho\leq\infty$, we set
$$\|u\|_{\tilde{L}^{\rho}_T(B^s_{p,r})}\overset{def}{=}\|2^{js}\|\Delta_ju\|_{L^{\rho}_T(L^p)}\|
_{l^r(\mathbb{Z})}.$$
We can then define the space $\tilde{L}^{\rho}_T(B^s_{p,r})$ as the set of tempered distributions
$u$ over $(0,T)\times \mathbb{R}^d$ such that $\|u\|_{\tilde{L}^{\rho}_T(B^s_{p,r})}\leq\infty$.
By the Minkowski inequality, we have
$$\|u\|_{\tilde{L}^{\rho}_T(B^s_{p,r})}\leq\|u\|_{L^{\rho}_T(B^s_{p,r})}~if~r\geq\rho$$
$$\|u\|_{L^{\rho}_T(B^s_{p,r})}\leq\|u\|_{\tilde{L}^{\rho}_T(B^s_{p,r})}~if~r\leq\rho.$$
The general principle is that all the properties of continuity for the product, composition, remainder, and
paraproduct remain true in those space.

Moreover when $s>0,~1\leq p\leq\infty,~1\leq \rho,\rho_1,\rho_2,\rho_3,\rho_4\leq\infty,$ and
$$\frac{1}{\rho}=\frac{1}{\rho_1}+\frac{1}{\rho_2}=\frac{1}{\rho_3}+\frac{1}{\rho_4},$$
we have
$$\|uv\|_{\tilde{L}^{\rho}_T(B^s_{p,r})}\leq C\left(\|u\|_{\tilde{L}^{\rho_1}_T(L^{\infty})}
\|v\|_{\tilde{L}^{\rho_2}_T(B^s_{p,r})}+\|v\|_{\tilde{L}^{\rho_3}_T(L^{\infty})}
\|u\|_{\tilde{L}^{\rho_4}_T(B^s_{p,r})}\right).$$
\end{lemm}

\begin{lemm}\label{h} \cite {Hs}
Les $s>1$, $u_0,h_0\in H^s$. Then there exist a positive time T, a unique solution $(u,h)$ of the Cauchy problem
(2.1) such that
$$u,h\in L^{\infty}([0,T],H^s),\nabla u\in L^2([0,T];H^s).$$
Furthermore, there exists a constant $c$ such that if $\|u_0\|_{H^s}+\|h_0\|_{H^s}\leq c$, then $T=\infty$.
\end{lemm}

\begin{rema}
All the proofs of Lemmas 2.6-2.16 can be found in \cite{H.J}.
\end{rema}

\section{The local well-posedness of Theorem \ref{t1}}
In order to study the local existence of solution, we define the function set $(u,h)\in\chi([0,T], s,p,r, E_1, E_2)$,  if $(u,h)\in\tilde{L}^{\infty}([0,T];B^s_{p,r})$, and
$$\|u\|_{\tilde{L}^{\infty}([0,T];B^s_{p,r})}\leq E_1,~~\|u\|_{ \tilde{L}^2([0,T];B^{s+1}_{p,r})}\leq E_1,~~\|h\|_{\tilde{L}^{\infty}([0,T];B^s_{p,r})}\leq E_2,$$
where $$E_1=8\nu^{-1}C_0\|u_0\|_{B^s_{p,r}},~E_2=4C_0\|h_0\|_{B^s_{p,r}}.$$

Next, we will prove Theorem \ref{t1} by the method of successive approximations.
Let us define the sequence $(u_n,h_n)$ by the following linear system:
\begin{align}
\left\{
\begin{array}{l}
(u_1,h_1)=S_2(u_0,h_0),\\[1ex]
\partial_tu_{n+1}+(u_n\cdot\nabla)u_{n+1}-\nu\Delta u_{n+1}=\frac{\nu}{1+h_n}\nabla h_n\nabla u_n+\nabla h_n,\\[1ex]
\partial_th_{n+1}+(u_n\cdot\nabla)h_{n+1}=-div\,u_n-h\,div\,u_n,\\[1ex]
(u_{n+1},h_{n+1})_{|t=0}=S_{n+2}(u_0,h_0).
\end{array}
\right.
\end{align}

Since $S_q$ are smooth operators, the initial date $S_{n+2}(u_0,h_0)$ are smooth functions. If $(u_n,h_n)\in \chi ([0,T],s,p,r,E_1,E_2)$ are smooth, then we have that for any $t\in[0,T]$,
$$\|h_n\|_{L^\infty} \leq C_{s,p} \|h_n\|_{B^s_{p,r}} \leq C_{s,p}E_2=4C_0C_{s,p}\|h_0\|_{B^s_{p,r}} \leq \frac{4C_0C_{s,p}}{8C_0C_{s,p}}=\frac{1}{2}.$$
Thus $\frac{\nu}{1+h_n}\nabla h_n\nabla u_n+\nabla h_n$ and $-div\,u_n-h\,div\,u_n$ are also smooth functions.
Note that the first equation in (3.1) is a transport diffusion equation for $u_{n+1}$, and the second equation is a transport
equation for $h_{n+1}$. Then the local existence of the smooth function for the Cauchy problem (3.1) is obvious.

 We split the proof of Theorem 1.1 into two steps: Estimation for big norms and Convergence for small norms.
 And for the sake of convenience, we suppose that $s<1+\frac{2}{p}$ and $p\neq\infty$($s\geq 1+\frac{2}{p}$ or $p=\infty$ is similar and easier).
\subsection{Estimation for big norms}
In this subsection, we want to prove the following proposition.
\begin{prop}\label{p1}
Suppose that $(u_0,h_0)\in B^s_{p,r}\times B^s_{p,r}$, $s>max\{1,\frac{2}{p}\}$, $1\leq p,r\leq\infty$, $\|h_0\|_{B^s_{p,r}}\leq \frac{1}{8CC_{s,p}}$, then
 there exists a positive time $T_1$, such that for any $n\in N$, $(u_n,h_n)\in\chi ([0,T_1],s,p,r,E_1,E_2)$.
\end{prop}
Proof:
Let $T(\geq T_1)$ satisfy
$$T\leq1,~e^{C_0^2E_1T}\leq2,~e^{2C_0C_{s,p}E_1T}\leq2,~(1+\nu T)^{\frac{3}{2}}\leq2.$$
 Then we prove the proposition by induction. Firstly let $(u_1,h_1)=S_2(u_0,h_0)$, thus we have
$$\|u_1\|_{\tilde{L}^{\infty}_{T_1}(B^s_{p,r})}\leq\|u_0\|_{B^S_{p,r}}\leq E_1,~
\|h_1\|_{\tilde{L}^{\infty}_{T_1}(B^s_{p,r})}\leq\|h_0\|_{B^S_{p,r}}\leq E_2,$$
$$\|u_1\|_{\tilde{L}^2_{T_1}(B^{s+1}_{p,r})}\leq T_1^{\frac{1}{2}}\|S_2u_0\|_{B^{s+1}_{p,r}}
\leq4\|S_2u_0\|_{B^s_{p,r}}\leq E_1$$

If $$\|u_n\|_{\tilde{L}^{\infty}_{T_1}(B^s_{p,r})}\leq E_1,~\|u_n\|_{\tilde{L}^2_{T_1}(B^{s+1}_{p,r})}\leq E_1,~
\|h_n\|_{\tilde{L}^{\infty}_{T_1}(B^s_{p,r})}\leq E_2,$$

then for $h_{n+1}$, in the view of Lemmas \ref{l12} and \ref{l13}, for all $t\leq T_1$, we have
\begin{align}
\begin{array}{l}
\|h_{n+1}\|_{\tilde{L}_t^{\infty}(B^s_{p,r})}\leq\left(\|S_{n+2}h_0\|_{B^s_{p,r}}+\|div\,u_n\|_{\tilde{L}_t^1(B^s_{p,r})}
+\|h_ndiv\,u_n\|_{\tilde{L}_t^1(B^s_{p,r})}\right)\\
exp\big(C_0\int_0^t\|\nabla u_n(t')\|_{L^\infty\cap B^{\frac{2}{p}}_{p,\infty}}dt'\big)\\[1ex]
\leq2\big(\frac{E_2}{4}+t^{\frac{1}{2}}\|div\,u_n\|_{\tilde{L}_t^2(B^s_{p,r})}+t^{\frac{1}{2}}
\|h_n\|_{\tilde{L}_t^\infty(B^s_{p,r})}\|div\,u_n\|_{L_t^2(L^{\infty})}\\
+t^{\frac{1}{2}}\|div\,u_n\|_{\tilde{L}_t^2(B^s_{p,r})}\|h_n\|_{L_t^\infty(L^{\infty})}\big)\\[1ex]
\leq\frac{E_2}{2}+Ct^{\frac{1}{2}}\|u_n\|_{\tilde{L}^2_t(B^{s+1}_{p,r})}+Ct^{\frac{1}{2}}
\|h_n\|_{\tilde{L}_t^\infty(B^s_{p,r})}\|u_n\|_{\tilde{L}^2_t(B^{s+1}_{p,r})}\\[1ex]
\leq\frac{E_2}{2}+C(1+E_2)t^{\frac{1}{2}}\|u_n\|_{\tilde{L}^2_t(B^{s+1}_{p,r})},
\end{array}
\end{align}
here we use the fact that
$$exp\big(C_0\int_0^t\|\nabla u_n(t')\|_{L^\infty\cap B^{\frac{2}{p}}_{p,\infty}}dt'\big)$$
$$\leq exp\big(2C_0C_{s,p}t^{\frac{1}{2}}\|u_n\|_{\tilde{L}^2_T(B^{s+1}_{p,r})}$$
$$\leq 2.$$

Now, we estimate $\|u_n\|_{\tilde{L}^\infty_t(B^s_{p,r})}$. By Lemmas \ref{l6} and \ref{l14},  we get
\begin{align}
\begin{array}{l}
\|u_{n+1}\|_{\tilde{L}^\infty_t(B^s_{p,r})}\leq C_0exp\big(C_0
\int_0^t\|\nabla u_n\|_{L^\infty\cap B^{\frac{2}{p}}_{p,r}}dt'\big)\Big(\|S_{n+2}u_0\|_{B^s_{p,r}}\\[2ex]
+\nu^{-\frac{1}{2}}(1+\nu t)^{\frac{1}{2}}\big(\|\nu\frac{\nabla h_n\nabla u_n}{1+h_n}\|
_{\tilde{L}^2_t(B^{s-1}_{p,r})}+\|\nabla h_n\|_{\tilde{L}^2_t(B^{s-1}_{p,r})}\big)\Big)\\[2ex]
\leq 2C_0\|u_0\|_{B^s_{p,r}}+Ct^{\frac{1}{2}}\|h_n\|_{\tilde{L}^{\infty}_t(B^s_{p,r})}+C\|\nabla(\ln(1+h_n))\nabla u_n\|_{\tilde{L}^2_t(B^{s-1}_{p,r})},
\end{array}
\end{align}
then by Lemmas \ref{0}, \ref{l10} and \ref{l11}, we have
\begin{align}
\begin{array}{l}
\|\nabla(\ln(1+h_n))\nabla u_n\|_{\tilde{L}^2_t(B^{s-1}_{p,r})}\\[1ex]
\leq C\|\nabla(\ln(1+h_n))\|_{\tilde{L}^\infty_t(B^{s-1-\frac{2}{p}}_{\infty,\infty})}\|\nabla u_n\|_{\tilde{L}_t^2(B^{\frac{2}{p}}_{p,r})}+C\|\nabla u_n\|_{L^2_t(L^\infty)}\|\nabla(\ln(1+h_n))\|_{\tilde{L}^\infty_t(B^{s-1}_{p,r})}\\[1ex]
+C\|\nabla u_n\|_{\tilde{L}^2_t(B^{\frac{2}{p}}_{p,r})}\|\nabla(\ln(1+h_n))\|_{\tilde{L}^\infty_t(B^{s-1-\frac{2}{p}}_{\infty,\infty})}\\[1ex]
\leq C\|h_n\|_{\tilde{L}_t^\infty(B^s_{p,r})}\|\nabla u_n\|_{\tilde{L}^2_t(B^{\frac{2}{p}+\ep}_{p,r})}\\[1ex]
\leq C\|h_n\|_{\tilde{L}_t^\infty(B^s_{p,r})}\|u_n\|_{\tilde{L}^2_t(B^{\frac{1}{p}+\frac{s}{2}+1}_{p,r})}\\[1ex]
\leq Ct^{\frac{s}{2}-\frac{1}{p}}\|h_n\|_{\tilde{L}_t^\infty(B^s_{p,r})}
\|u_n\|_{\tilde{L}^{\frac{2}{\frac{1}{p}+1-\frac{s}{2}}}_t(B^{\frac{1}{p}+\frac{s}{2}+1}_{p,r})}\\[2ex]
\leq Ct^{\frac{s}{2}-\frac{1}{p}}\|h_n\|_{\tilde{L}_t^\infty(B^s_{p,r})}\|u_n\|^{\frac{1}{p}+1-\frac{s}{2}}
_{\tilde{L}^2_t(B^{s+1}_{p,r})}\|u_n\|^{\frac{s}{2}-\frac{1}{p}}_{\tilde{L}^\infty_t(B^s_{p,r})},
\end{array}
\end{align}
combining (3.3) and (3.4), we have

\begin{align}
\begin{array}{l}
\|u_{n+1}\|_{\tilde{L}^\infty_t(B^s_{p,r})}\leq 2C_0\|u_0\|_{B^s_{p,r}}+
Ct^{\frac{1}{2}}\|h_n\|_{\tilde{L}_t^\infty(B^s_{p,r})}\\[1ex]
+Ct^{\frac{s}{2}-\frac{1}{p}}\|h_n\|_{\tilde{L}_t^\infty(B^s_{p,r})}\|u_n\|^{\frac{1}{p}+1-\frac{s}{2}}
_{\tilde{L}^2_t(B^{s+1}_{p,r})}\|u_n\|^{\frac{s}{2}-\frac{1}{p}}_{\tilde{L}^\infty_t(B^s_{p,r})}.
\end{array}
\end{align}
Similarly as $\|u_{n+1}\|_{\tilde{L}^\infty_{B^s_{p,r}}}$, by Lemma \ref{l14}, we have
\begin{align}
\begin{array}{l}
\|u_{n+1}\|_{\tilde{L}^2_t(B^{s+1}_{p,r})}\leq 4\nu^{-\frac{1}{2}}C_0\|u_0\|_{B^s_{p,r}}+\\[1ex]
Ct^{\frac{1}{2}}\|h_n\|_{\tilde{L}_t^\infty(B^s_{p,r})}+Ct^{\frac{s}{2}-\frac{1}{p}}\|h_n\|_{\tilde{L}_t^\infty(B^s_{p,r})}\|u_n\|^{\frac{1}{p}+1-\frac{s}{2}}
_{\tilde{L}^2_t(B^{s+1}_{p,r})}\|u_n\|^{\frac{s}{2}-\frac{1}{p}}_{\tilde{L}^\infty_t(B^s_{p,r})}.
\end{array}
\end{align}

Thus let

$T'_1=min\{T,(4CE_1)^{-1},~\big(4C^2(1+E_2)(E_1+E_2+E_1E_2)\big)^{-2}E^2_2\},$ by (3.2), (3.5), (3.6),
we can get that, for any $t\leq T_1\leq T'_1$,
$$\|h_{n+1}\|_{\tilde{L}_t^{\infty}(B^s_{p,r})}\leq E_2,~
\|u_{n+1}\|_{\tilde{L}_t^{\infty}(B^s_{p,r})}\leq E_1,~
\|u_{n+1}\|_{\tilde{L}_t^2(B^{s+1}_{p,r})}\leq E_1.$$
This completes the proof of Proposition \ref{p1}.
\subsection{Convergence of small norms}
\begin{prop}\label{p2}
Suppose that $(u_0,h_0)\in B^s_{p,r}\times B^s_{p,r}$, $s>max\{\frac{2}{p},1\}$, $1\leq p,r\leq\infty$,  $\|h_0\|_{B^s_{p,r}}\leq \frac{1}{8CC_{s,p}}$, then
 there exists a positive time $T_2(\leq T_1)$, such that $(u_n,h_n)$ is a Cauchy sequence in $\chi ([0,T_2],s-1,p,r,E_1,E_2)$.
\end{prop}
Proof: From the equations in (3.1), we have
\begin{align}
\left\{
\begin{array}{l}
\partial_t(u_{n+1}-u_n)+(u_n\cdot\nabla)(u_{n+1}-u_n)-\nu\Delta(u_{n+1}-u_n)=\sum_{j=1}^5F_j\\[1ex]
\partial_t((h_{n+1}-h_n)+(u_n\cdot\nabla)(h_{n+1}-h_n)=\sum_{j=1}^4J_j,\\[1ex]
(u_{n+1}-u_n,h_{n+1}-h_n)_{|t=0}=\Delta_{n+1}(u_0,h_0),
\end{array}
\right.
\end{align}

where
\begin{align}
\begin{array}{l}
\sum_{j=1}^5F_j=(u_n-u_{n-1})\cdot\nabla u_n+\nabla(h_n-h_{n-1})+\frac{\nu}{1+h_n}\nabla h_n\nabla(u_n-u_{n-1})\\[1ex]
+\frac{\nu}{1+h_n}\nabla u_{n-1}\nabla(h_n-h_{n-1})+\nu
(\frac{1}{1+h_n}-\frac{1}{1+h_{n-1}})\nabla h_{n-1}\nabla u_{n-1},\\[2ex]
\sum_{j=1}^4J_j=(u_n-u_{n-1})\cdot\nabla h_n+div\,(u_n-u_{n-1})+h_n\,div\,(u_n-u_{n-1})\\[1ex]
+(h_n-h_{n-1})\,div\,u_{n-1}.
\end{array}
\end{align}

Then we estimate the Besov norm of $u_{n+1}-u_n$ and $h_{n+1}-h_n$.
For any $t\leq T_2\leq T_1$, by Lemma \ref{l14}, we have
\begin{align}
\begin{array}{l}
\|u_{n+1}-u_n\|_{\tilde{L}^\infty_t(B^{s-1}_{p,r})}\leq C_0exp\big(C_0\int_0^t\|\nabla u_n\|_{B^{\frac{2}{p}}_{p,\infty}\cap L^\infty}dt'\big)\times\\[1ex]
\left(\|S_{n+2}u_0-S_{n+1}u_0\|_{B^{s-1}_{p,r}}+(\frac{1+\nu t}{\nu})^{\frac{1}{2}}
\|\sum_{j=1}^5F_j\|_{\tilde{L}^2_t(B^{s-2}_{p,r})}\right)\\[1ex]
\leq 2C_0\|\Delta_{n+1}u_0\|_{B^{s-1}_{p,r}}+C
\|\sum_{j=1}^5F_j\|_{\tilde{L}^2_t(B^{s-2}_{p,r})},
\end{array}
\end{align}
here
\begin{align}
\begin{array}{l}
\|\sum_{j=1}^5F_j\|_{\tilde{L}^2_t(B^{s-2}_{p,r})}\leq
\|(u_n-u_{n-1})\cdot\nabla u_n\|_{\tilde{L}^2_t(B^{s-2}_{p,r})}\\[1ex]
+\|\nabla(h_n-h_{n-1})\|_{\tilde{L}^2_t(B^{s-2}_{p,r})}
+\|\nu\frac{\nabla h_n}{1+h_n}\nabla(u_n-u_{n-1})\|_{\tilde{L}^2_t(B^{s-2}_{p,r})}\\[1ex]
+\|\nu\frac{\nabla(h_n-h_{n-1})}{1+h_n}\nabla u_{n-1}\|_{\tilde{L}^2_t(B^{s-2}_{p,r})}
+\|\nu\frac{h_n-h_{n-1}}{(1+h_n)(1+h_{n-1})}\nabla h_{n-1}\nabla u_{n-1}\|_{\tilde{L}^2_t(B^{s-2}_{p,r})}\\[1ex]
=I_1+I_2+I_3+I_4+I_5,
\end{array}
\end{align}
Next, we deal with $I_j,~j=1,2,3,4,5$ term by term.

By Lemmas \ref{l5}-\ref{l6}, Lemma \ref{l10} and Lemma \ref{l11}, we have
\begin{align}
\begin{array}{l}
I_1=\|(u_n-u_{n-1})\cdot\nabla u_n\|_{\tilde{L}^2_t(B^{s-2}_{p,r})}\\[1ex]
\leq C\|u_n-u_{n-1}\|_{\tilde{L}^{\infty}_t(B^{s-1-\frac{2}{p}}_{\infty,\infty})}\|\nabla u_n\|_{\tilde{L}_t^2(B^{\frac{2}{p}-1}_{p,r})}+
C\|\nabla u_n\|_{\tilde{L}^2_t(B^{-1}_{\infty,\infty})}\|u_n-u_{n-1}\|_{L_t^{\infty}(B^{s-1}_{p,r})}\\[1ex]
+C\|R\big(\nabla u_n,(u_n-u_{n-1})\big)\|_{\tilde{L}^2_t(B^{s-2}_{p,r})}.
\end{array}
\end{align}
If $p\leq 2$, we have
\begin{align}
\begin{array}{l}
\|R\big(\nabla u_n,(u_n-u_{n-1})\big)\|_{\tilde{L}^2_t(B^{s-2}_{p,r})}\\[1ex]
\leq C\|R\big(\nabla u_n,(u_n-u_{n-1})\big)\|_{\tilde{L}^2_t(B^{s-\frac{2}{p}}_{1,r})}\\[1ex]
\leq C\|u_n-u_{n-1}\|_{\tilde{L}_t^\infty(B^{s-1}_{p,r})}\|\nabla u_n\|_{\tilde{L}^2_t(B^{1-\frac{2}{p}}_{p',\infty})}\\[1ex]
\leq C\|u_n-u_{n-1}\|_{\tilde{L}_t^\infty(B^{s-1}_{p,r})}\|\nabla u_n\|_{\tilde{L}^2_t(B^{1-\frac{2}{p}+\frac{2}{p}-\frac{2}{p'}}_{p,\infty})}\\[1ex]
\leq C\|u_n-u_{n-1}\|_{\tilde{L}_t^\infty(B^{s-1}_{p,r})}\|u_n\|_{\tilde{L}^2_t(B^{\frac{2}{p}}_{p,r})},
\end{array}
\end{align}
where $\frac{1}{p'}=1-\frac{1}{p}$.

If $p>2$ and $s-2+\frac{2}{p}\leq0$, we have
\begin{align}
\begin{array}{l}
\|R\big(\nabla u_n,(u_n-u_{n-1})\big)\|_{\tilde{L}^2_t(B^{s-2}_{p,r})}\\[1ex]
\leq C\|R\big(\nabla u_n,(u_n-u_{n-1})\big)\|_{\tilde{L}^2_t(B^{s-2+\frac{2}{p}}_{\frac{p}{2},r})}\\[1ex]
\leq C\|R\big(\nabla u_n,(u_n-u_{n-1})\big)\|_{\tilde{L}^2_t(B^\ep_{\frac{p}{2},r})}\\[1ex]
\leq C\|u_n-u_{n-1}\|_{\tilde{L}_t^\infty(B^{s-1}_{p,r})}\|\nabla u_n\|_{\tilde{L}^2_t(B^{1-s+\ep}_{p,\infty})}\\[1ex]
\leq C\|u_n-u_{n-1}\|_{\tilde{L}_t^\infty(B^{s-1}_{p,r})}\| u_n\|_{\tilde{L}^2_t(B^{2-s+\ep}_{p,\infty})}\\[1ex]
\leq C\|u_n-u_{n-1}\|_{\tilde{L}_t^\infty(B^{s-1}_{p,r})}\| u_n\|_{\tilde{L}^2_t(B^s_{p,\infty})},
\end{array}
\end{align}
here hereafter, $\ep$ is a positive small real constant and it does not influence the direct of inequalities, for example if $s<2$, then $s+\ep<2$ as well.
And here we also use the fact $2-s+\ep<s$ by $s>1$.

If $s-2+\frac{2}{p}>0$, we have the similar estimate as (3.12), then combining (3.11)-(3.13), we can obtain
\begin{align}
\begin{array}{l}
I_1=\|(u_n-u_{n-1})\cdot\nabla u_n\|_{\tilde{L}^2_t(B^{s-2}_{p,r})}
\leq C\|u_n\|_{\tilde{L}^2_t(B^s_{p,r})}\|u_n-u_{n-1}\|_{\tilde{L}^{\infty}_t(B^{s-1}_{p,r})}\\[1ex]
\leq C E_1 t^{\frac{1}{2}}\|u_n-u_{n-1}\|_{\tilde{L}^{\infty}_t(B^{s-1}_{p,r})}.
\end{array}
\end{align}
From Lemma \ref{l6}, it's easy to see that
\begin{align}
\begin{array}{l}
I_2\leq t^{\frac{1}{2}}\|\nabla(h_n-h_{n-1})\|_{\tilde{L}^{\infty}_t(B^{s-2}_{p,r})}\leq
Ct^{\frac{1}{2}}\|h_n-h_{n-1}\|_{\tilde{L}^{\infty}_t(B^{s-1}_{p,r})}.
\end{array}
\end{align}
In view of Lemmas \ref{l5}-\ref{l6}, Lemmas \ref{l10}-\ref{l11}, we get
\begin{align}
\begin{array}{l}
I_3=\|\nu\frac{\nabla h_n}{1+h_n}\nabla(u_n-u_{n-1})\|_{\tilde{L}^2_t(B^{s-2}_{p,r})}\\[1ex]
=\|\nu\nabla(\ln(1+h_n))\nabla(u_n-u_{n-1})\|_{\tilde{L}^{\infty}_t(B^{s-2}_{p,r})}\\[1ex]
\leq C\|\nabla(\ln(1+h_n))\|_{\tilde{L}^\infty_t(B^{s-1-\frac{2}{p}}_{\infty,\infty})}\|\nabla(u_n-u_{n-1})\|
_{\tilde{L}^2_t(B^{\frac{2}{p}-1}_{p,r})}\\[1ex]
+C\|\nabla(u_n-u_{n-1})\|
_{\tilde{L}^2_t(B^{-1}_{\infty,\infty})}\|\nabla(\ln(1+h_n))\|_{\tilde{L}^\infty_t(B^{s-1}_{p,r})}\\[1ex]
+C\|R\big(\nabla(\ln(1+h_n)),\nabla(u_n-u_{n-1})\|_{\tilde{L}^2_t(B^{s-2}_{p,r})}.
\end{array}
\end{align}
Similarly as $I_1$, we have
\begin{align}
\begin{array}{l}
\|R\big(\nabla(\ln(1+h_n)),\nabla(u_n-u_{n-1})\|_{\tilde{L}^2_t(B^{s-2}_{p,r})}\\[1ex]
\leq C\|h_n\|_{\tilde{L}^\infty_t(B^s_{p,r})}\big(\|u_n-u_{n-1}\|_{\tilde{L}^2_t(B^{\frac{2}{p}}_{p,r})}+
\|u_n-u_{n-1}\|_{\tilde{L}^2_t(B^{2-s+\ep}_{p,r})}\big)\\[1ex]
\leq C\|h_n\|_{\tilde{L}^\infty_t(B^s_{p,r})}\big(\|u_n-u_{n-1}\|_{\tilde{L}^2_t(B^{\frac{2}{p}}_{p,r})}+
\|u_n-u_{n-1}\|_{\tilde{L}^2_t(B^1_{p,r})}\big).
\end{array}
\end{align}
Combining (3.16) and (3.17), we have
\begin{align}
\begin{array}{l}
I_3=\|\nu\frac{\nabla h_n}{1+h_n}\nabla(u_n-u_{n-1})\|_{\tilde{L}^2_t(B^{s-2}_{p,r})}\\[1ex]
\leq C\|h_n\|_{\tilde{L}^\infty_t(B^s_{p,r})}\big(\|u_n-u_{n-1}\|_{\tilde{L}^2_t(B^{\frac{2}{p}}_{p,r})}+
\|u_n-u_{n-1}\|_{\tilde{L}^2_t(B^1_{p,r})}\big).
\end{array}
\end{align}
By $s-1<\frac{2}{p}<s$ and in the view of Lemma \ref{0}, we have
\begin{align}
\begin{array}{l}
\|u_n-u_{n-1}\|_{\tilde{L}^2_t(B^{\frac{2}{p}}_{p,r})}\\[1ex]
\leq t^{s-\frac{2}{p}}\|u_n-u_{n-1}\|_{\tilde{L}^{\frac{2}{1+\frac{2}{p}-s}}_t(B^{\frac{2}{p}}_{p,r})}\\[1ex]
\leq t^{s-\frac{2}{p}}\|u_n-u_{n-1}\|^{1+\frac{2}{p}-s}_{\tilde{L}^2_t(B^s_{p,r})}
\|u_n-u_{n-1}\|^{s-\frac{2}{p}}_{\tilde{L}^\infty_t(B^{s-1}_{p,r})}\\[1ex]
\leq t^{s-\frac{2}{p}}\big(\|u_n-u_{n-1}\|_{\tilde{L}^2_t(B^s_{p,r})}+
\|u_n-u_{n-1}\|_{\tilde{L}^\infty_t(B^{s-1}_{p,r})}\big).
\end{array}
\end{align}
As regards $\|u_n-u_{n-1}\|_{\tilde{L}^2_t(B^1_{p,r})}$, if $s\ge2$, we have
\begin{align}
\begin{array}{l}
\|u_n-u_{n-1}\|_{\tilde{L}^2_t(B^1_{p,r})}\le t^{\frac{1}{2}}\|u_n-u_{n-1}\|_{\tilde{L}^\infty_t(B^{s-1}_{p,r})},
\end{array}
\end{align}
if $s<2$, it means $s-1<1<s$, by Lemma \ref{0}, we have
\begin{align}
\begin{array}{l}
\|u_n-u_{n-1}\|_{\tilde{L}^2_t(B^1_{p,r})}\\[1ex]
\leq t^{\frac{s-1}{2}}\|u_n-u_{n-1}\|_{\tilde{L}^{\frac{2}{2-s}}_t(B^1_{p,r})}\\[1ex]
\leq t^{\frac{s-1}{2}}\|u_n-u_{n-1}\|^{2-s}_{\tilde{L}^2_t(B^s_{p,r})}
\|u_n-u_{n-1}\|^{s-1}_{\tilde{L}^\infty_t(B^{s-1}_{p,r})}\\[1ex]
\leq t^{\frac{s-1}{2}}\big(\|u_n-u_{n-1}\|_{\tilde{L}^2_t(B^s_{p,r})}+
\|u_n-u_{n-1}\|_{\tilde{L}^2_t(B^{s-1}_{p,r})}\big).
\end{array}
\end{align}
Thus we have
\begin{align}
\begin{array}{l}
I_3=\|\nu\frac{\nabla h_n}{1+h_n}\nabla(u_n-u_{n-1})\|_{\tilde{L}^2_t(B^{s-2}_{p,r})}\\[1ex]
\leq CE_2\big(t^\frac{1}{2}+t^{s-\frac{2}{p}}+t^{\frac{s-1}{2}}\big)
\big(\|u_n-u_{n-1}\|_{\tilde{L}^2_t(B^s_{p,r})}+
\|u_n-u_{n-1}\|_{\tilde{L}^2_t(B^{s-1}_{p,r})}\big).
\end{array}
\end{align}

Then we deal with $I_4$ by the similar method.
\begin{align}
\begin{array}{l}
I_4=\|\nu\frac{\nabla(h_n-h_{n-1})}{1+h_n}\nabla u_{n-1}\|_{\tilde{L}^2_t(B^{s-2}_{p,r})}\\[1ex]
=\nu\|(1-\frac{h_n}{1+h_n})\nabla(h_n-h_{n-1})\nabla u_{n-1}\|_{\tilde{L}^2_t(B^{s-2}_{p,r})}\\[1ex]
\leq\nu\|\nabla(h_n-h_{n-1})\nabla u_{n-1}\|_{\tilde{L}^2_t(B^{s-2}_{p,r})}+\nu\|\frac{h_n}{1+h_n}\nabla(h_n-h_{n-1})\nabla u_{n-1}\|_{\tilde{L}^2_t(B^{s-2}_{p,r})}\\[1ex]
\leq I_{41}+I_{42}.
\end{array}
\end{align}

Similar to the argument in the proof of $I_3$, we obtain
\begin{align}
\begin{array}{l}
I_{41}\leq C\|\nabla(h_n-h_{n-1})\|_{\tilde{L}^\infty_t(B^{s-2}_{p,r})}\big(\|\nabla u_{n-1}\|_{\tilde{L}^2_t(B^{\frac{2}{p}+\ep}_{p,r})}+|\nabla u_{n-1}\|_{\tilde{L}^2_t(B^{2-s+\ep}_{p,r})}\big)\\[1ex]
\leq  C\|h_n-h_{n-1}\|_{\tilde{L}^\infty_t(B^{s-1}_{p,r})}\big(\| u_{n-1}\|_{\tilde{L}^2_t(B^{\frac{1}{p}+\frac{s}{2}+1}_{p,r})}+\| u_{n-1}\|_{\tilde{L}^2_t(B^2_{p,r})}\big).
\end{array}
\end{align}
Following the procedure of (3.19)-(3.21) respectively, we have
\begin{align}
\begin{array}{l}
\| u_{n-1}\|_{\tilde{L}^2_t(B^{\frac{1}{p}+\frac{s}{2}+1}_{p,r})}\leq Ct^{\frac{s}{2}-\frac{1}{p}}\big(\|u_{n-1}\|_{\tilde{L}^\infty_t(B^s_{p,r})}+\|u_{n-1}\|_{\tilde{L}^2_t(B^{s+1}_{p,r})}\big),
\end{array}
\end{align}
and in the case $s\ge2$
\begin{align}
\begin{array}{l}
\| u_{n-1}\|_{\tilde{L}^2_t(B^2_{p,r})}\leq Ct^{\frac{1}{2}}\|u_{n-1}\|_{\tilde{L}^\infty_t(B^s_{p,r})},
\end{array}
\end{align}
in the case $s<2$
\begin{align}
\begin{array}{l}
\| u_{n-1}\|_{\tilde{L}^2_t(B^2_{p,r})}\leq Ct^{\frac{s-1}{2}}\big(\|u_{n-1}\|_{\tilde{L}^\infty_t(B^s_{p,r})}+\|u_{n-1}\|_{\tilde{L}^2_t(B^{s+1}_{p,r})}\big).
\end{array}
\end{align}
Thus we have
\begin{align}
\begin{array}{l}
I_{41}
\leq C\big(t^{\frac{s}{2}-\frac{1}{p}}+t^{\frac{1}{2}}+t^{\frac{s-1}{2}}\big)  \|h_n-h_{n-1}\|_{\tilde{L}^\infty_t(B^{s-1}_{p,r})}\big(\|u_{n-1}\|_{\tilde{L}^\infty_t(B^s_{p,r})}+\|u_{n-1}\|_{\tilde{L}^2_t(B^{s+1}_{p,r})}\big)\\[1ex]
\leq CE_1\big(t^{\frac{s}{2}-\frac{1}{p}}+t^{\frac{1}{2}}+t^{\frac{s-1}{2}}\big)  \|h_n-h_{n-1}\|_{\tilde{L}^\infty_t(B^{s-1}_{p,r})}.
\end{array}
\end{align}
By Lemmas \ref{l10}-\ref{l11},
we have
\begin{align}
\begin{array}{l}
\|\frac{h_n}{1+h_n}\nabla u_{n-1}\|_{\tilde{L}^2_t(B^{\frac{2}{p}+\ep}_{p,r})}
\leq C\|h_n\|_{\tilde{L}^\infty_t(B^s_{p,r})}\|\nabla u_{n-1}\|_{\tilde{L}^2_t(B^{\frac{2}{p}+\ep}_{p,r})},
\end{array}
\end{align}
and
\begin{align}
\begin{array}{l}
\|\frac{h_n}{1+h_n}\nabla u_{n-1}\|_{\tilde{L}^2_t(B^{2-s+\ep}_{p,r})}
\leq C\|h_n\|_{\tilde{L}^\infty_t(B^s_{p,r})}\big(\|\nabla u_{n-1}\|_{\tilde{L}^2_t(B^{\frac{2}{p}+\ep}_{p,r})}+\|\nabla u_{n-1}\|_{\tilde{L}^2_t(B^{2-s+\ep}_{p,r})}\big).
\end{array}
\end{align}
Thus we have
\begin{align}
\begin{array}{l}
I_{42}
\leq CE_1E_2\big(t^{\frac{s}{2}-\frac{1}{p}}+t^{\frac{1}{2}}+t^{\frac{s-1}{2}}\big)  \|h_n-h_{n-1}\|_{\tilde{L}^\infty_t(B^{s-1}_{p,r})}.
\end{array}
\end{align}
Combining (3.28) and (3.31), we have
\begin{align}
\begin{array}{l}
I_4
\leq CE_1(E_2+1)\big(t^{\frac{s}{2}-\frac{1}{p}}+t^{\frac{1}{2}}+t^{\frac{s-1}{2}}\big)  \|h_n-h_{n-1}\|_{\tilde{L}^\infty_t(B^{s-1}_{p,r})}.
\end{array}
\end{align}

Similarly as $I_4$, we have
\begin{align}
\begin{array}{l}
I_5
\leq CE_1E_2(E_2+1)\big(t^{\frac{s}{2}-\frac{1}{p}}+t^{\frac{1}{2}}+t^{\frac{s-1}{2}}\big)  \|h_n-h_{n-1}\|_{\tilde{L}^\infty_t(B^{s-1}_{p,r})}.
\end{array}
\end{align}

We also have
\begin{align}
\|\Delta_{n+1}u_0\|_{B^{s-1}_{p,r}}\leq2^{-(n+1)}\|\Delta_{n+1}u_0\|_{B^s_{p,r}}
\leq2^{-(n+1)}\|u_0\|_{B^s_{p,r}}.
\end{align}
Combining (3.9)-(3.34), we have
\begin{align}
\begin{array}{l}
\|u_{n+1}-u_n\|_{\tilde{L}^\infty_t(B^{s-1}_{p,r})}
\leq 2C_02^{-(n+1)}\|u_0\|_{B^s_{p,r}}\\[1ex]
+C(1+E_1+E_1E_2+E_1E_2^2)\big(t^{\frac{s}{2}-\frac{1}{p}}+t^{\frac{1}{2}}+t^{\frac{s-1}{2}}\big)  \|h_n-h_{n-1}\|_{\tilde{L}^\infty_t(B^{s-1}_{p,r})}\\[1ex]
C(E_1+E_2)\big(t^\frac{1}{2}+t^{s-\frac{2}{p}}+t^{\frac{s-1}{2}}\big)
\big(\|u_n-u_{n-1}\|_{\tilde{L}^2_t(B^s_{p,r})}+
\|u_n-u_{n-1}\|_{\tilde{L}^2_t(B^{s-1}_{p,r})}\big).
\end{array}
\end{align}
For $\|u_{n+1}-u_n\|_{\tilde{L}^2_t(B^s_{p,r})}$, we have the similar estimate
\begin{align}
\begin{array}{l}
\|u_{n+1}-u_n\|_{\tilde{L}^2_t(B^s_{p,r})}\leq C_0\nu^{-\frac{1}{2}}exp\big(C_0(1+\nu t)^{\frac{1}{2}}\int_0^t\|\nabla u_n\|_{B^{\frac{2}{p}}_{p,\infty}\cap L^\infty}dt'\big)\times\\[1ex]
\left((1+\nu t)^{\frac{1}{2}}\|S_{n+2}u_0-S_{n+1}u_0\|_{B^{s-1}_{p,r}}+\frac{1+\nu t}{\nu}
\|\sum_{j=1}^5F_j\|_{\tilde{L}^2_t(B^{s-1}_{p,r})}\right)\\[1ex]
\leq 4\nu^{-1}C_02^{-(n+1)}\|u_0\|_{B^s_{p,r}}\\[1ex]
+C(1+E_1+E_1E_2+E_1E_2^2)\big(t^{\frac{s}{2}-\frac{1}{p}}+t^{\frac{1}{2}}+t^{\frac{s-1}{2}}\big)  \|h_n-h_{n-1}\|_{\tilde{L}^\infty_t(B^{s-1}_{p,r})}\\[1ex]
C(E_1+E_2)\big(t^\frac{1}{2}+t^{s-\frac{2}{p}}+t^{\frac{s-1}{2}}\big)
\big(\|u_n-u_{n-1}\|_{\tilde{L}^2_t(B^s_{p,r})}+
\|u_n-u_{n-1}\|_{\tilde{L}^2_t(B^{s-1}_{p,r})}\big).
\end{array}
\end{align}

For $h_{n+1}-h_n$, we have
\begin{align}
\begin{array}{l}
\|h_{n+1}-h_n\|_{\tilde{L}^{\infty}_t(B^{s-1}_{p,r})}\leq exp(C_0\int_0^t\|\nabla u_n\|_{B^{\frac{2}{p}}_{p,\infty}\cap L^\infty}dt')\times \\[1ex]
\big(\|\Delta_{n+1}h_0\|_{B^{s-1}_{p,r}}
+\|\sum_{j=1}^4J_j\|_{\tilde{L}^1_t(B^{s-1}_{p,r})}\big)\\[1ex]
\leq2\big(\|\Delta_{n+1}h_0\|_{B^{s-1}_{p,r}}
+t^{\frac{1}{2}}
\|\sum_{j=1}^4J_j\|_{\tilde{L}^2_t(B^{s-1}_{p,r})}\big).
\end{array}
\end{align}

From Lemmas \ref{l5}-\ref{l6}, Lemmas \ref{l10}-\ref{l11}, we have
\begin{align}
\begin{array}{l}
\|J_1\|_{\tilde{L}_t^2(B^{s-1}_{p,r})}=\|(u_n-u_{n-1})\cdot\nabla h_n\|_{\tilde{L}_t^{\infty}(B^{s-1}_{p,r})}\\[1ex]
\leq \|u_n-u_{n-1}\|_{L_t^2(L^\infty)}\|\nabla h_n\|_{\tilde{L}_t^{\infty}(B^{s-1}_{p,r})}
+\|\nabla h_n\|_{\tilde{L}_t^{\infty}(B^{s-1-\frac{2}{p}}_{\infty,\infty})}\|u_n-u_{n-1}\|_{\tilde{L}_t^2(B^{\frac{2}{p}}_{p,r})}\\[1ex]
\leq C\|h_n\|_{\tilde{L}_t^{\infty}(B^s_{p,r})}\|u_n-u_{n-1}\|_{\tilde{L}_t^2(B^s_{p,r})}\\
\leq CE_2\|u_n-u_{n-1}\|_{\tilde{L}_t^2(B^s_{p,r})}.
\end{array}
\end{align}
In view of Lemma \ref{l6}, we get
\begin{align}
\begin{array}{l}
\|J_2\|_{\tilde{L}_t^2(B^{s-1}_{p,r})}=\|div\,(u_n-u_{n-1})\|_{\tilde{L}_t^2(B^{s-1}_{p,r})}\\[1ex]
\leq C\|u_n-u_{n-1}\|_{\tilde{L}_t^2(B^s_{p,r})}.
\end{array}
\end{align}
By Lemmas \ref{l5}-\ref{l6}, Lemmas \ref{l10}-\ref{l11}, we obtain
\begin{align}
\begin{array}{l}
\|J_3\|_{\tilde{L}_t^2(B^{s-1}_{p,r})}=\|h_n\,div\,(u_n-u_{n-1})\|_{\tilde{L}_t^2(B^{s-1}_{p,r})}\\[1ex]
\leq\|h_n\|_{L^{\infty}_t(L^{\infty})}\|div\,(u_n-u_{n-1})\|_{\tilde{L}_t^2(B^{s-1}_{p,r})}
+\|div\,(u_n-u_{n-1})\|_{L^2_t(L^{\infty})}\|h_n\|_{\tilde{L}_t^{\infty}(B^{s-1}_{p,r})}\\[1ex]
\leq C\|h_n\|_{\tilde{L}_t^{\infty}(B^{s-1}_{p,r})}\|u_n-u_{n-1}\|_{\tilde{L}_t^2(B^s_{p,r})}\\
\leq CE_2\|u_n-u_{n-1}\|_{\tilde{L}_t^2(B^s_{p,r})},
\end{array}
\end{align}
Similarly to $J_3$, we have
\begin{align}
\begin{array}{l}
\|J_4\|_{\tilde{L}_t^2(B^{s-1}_{p,r})}=\|(h_n-h_{n-1})\,div\,u_{n-1}\|_{\tilde{L}_t^{\infty}(B^{s-1}_{p,r})}\\[1ex]
\leq\|h_n-h_{n-1}\|_{\tilde{L}^{\infty}_t(B^{s-1-\frac{2}{p}}_{\infty,\infty}}\|div\,u_{n-1}\|_{\tilde{L}_t^2(B^{\frac{2}{p}}_{p,r})}
+\|div\,u_{n-1}\|_{L^2_t(L^{\infty})}\|h_n-h_{n-1}\|_{\tilde{L}_t^{\infty}(B^{s-1}_{p,r})}\\[1ex]
\leq C\|h_n-h_{n-1}\|_{\tilde{L}_t^{\infty}(B^{s-1}_{p,r})}\|u_{n-1}\|_{\tilde{L}_t^2(B^{s+1}_{p,r})}\\
\leq CE_1\|h_n-h_{n-1}\|_{\tilde{L}_t^{\infty}(B^{s-1}_{p,r})}.
\end{array}
\end{align}
We also have
\begin{align}
\|\Delta_{n+1}h_0\|_{B^{s-1}_{p,r}}\leq2^{-(n+1)}\|\Delta_{n+1}h_0\|_{B^s_{p,r}}
\leq2^{-(n+1)}\|h_0\|_{B^s_{p,r}}.
\end{align}
Thus we have
\begin{align}
\begin{array}{l}
\|h_{n+1}-h_n\|_{\tilde{L}^{\infty}_t(B^{s-1}_{p,r})}\leq 2^{-n}\|h_0\|_{B^s_{p,r}}\\
+C(1+E_2)t^\frac{1}{2}\|u_n-u_{n-1}\|_{\tilde{L}^2_t(B^s_{p,r})}
+CE_1 t^{\frac{1}{2}}\|h_n-h_{n-1}\|_{\tilde{L}^\infty_t(B^{s-1}_{p,r})}.
\end{array}
\end{align}

Choose a suitable $T_2(\leq T_1)$ such that:
\begin{align}
\left\{
\begin{array}{l}
C(1+E_1+E_1E_2+E_1E_2^2)\big(T_2^{\frac{s}{2}-\frac{1}{p}}+T_2^{\frac{1}{2}}+T_2^{\frac{s-1}{2}}\big)\leq\frac{1}{12},\\[1ex]
C(E_1+E_2)\big(T_2^\frac{1}{2}+T_2^{s-\frac{2}{p}}+T_2^{\frac{s-1}{2}}\big)\leq\frac{1}{12},\\[1ex]
C(1+E_2)T_2^\frac{1}{2}\leq\frac{1}{12},~
+CE_1 T_2^{\frac{1}{2}}\leq\frac{1}{12}.
\end{array}
\right.
\end{align}
Thus, for any $t\leq T_2$, we can obtain
\begin{align}
\begin{array}{l}
\|u_{n+1}-u_n\|_{\tilde{L}^{\infty}_t(B^{s-1}_{p,r})}\leq\frac{1}{4}2^{-n}E_1\\[1ex]
+\frac{1}{12}\|u_n-u_{n-1}\|_{\tilde{L}^{\infty}_t(B^{s-1}_{p,r})}+
\frac{1}{12}\|u_n-u_{n-1}\|_{\tilde{L}^2_t(B^s_{p,r})}+
\frac{1}{12}\|h_n-h_{n-1}\|_{\tilde{L}^{\infty}_t(B^{s-1}_{p,r})},
\end{array}
\end{align}
\begin{align}
\begin{array}{l}
\|u_{n+1}-u_n\|_{\tilde{L}^2_t(B^s_{p,r})}\leq\frac{1}{4}2^{-n}E_1\\[1ex]
+\frac{1}{12}\|u_n-u_{n-1}\|_{\tilde{L}^{\infty}_t(B^{s-1}_{p,r})}+
\frac{1}{12}\|u_n-u_{n-1}\|_{\tilde{L}^2_t(B^s_{p,r})}+
\frac{1}{12}\|h_n-h_{n-1}\|_{\tilde{L}^{\infty}_t(B^{s-1}_{p,r})},
\end{array}
\end{align}
and
\begin{align}
\begin{array}{l}
\|h_{n+1}-h_n\|_{\tilde{L}^{\infty}_t(B^{s-1}_{p,r})}\leq\frac{1}{4}2^{-n}E_2\\[1ex]
+\frac{1}{12}\|u_n-u_{n-1}\|_{\tilde{L}^2_t(B^s_{p,r})}+
\frac{1}{12}\|h_n-h_{n-1}\|_{\tilde{L}^{\infty}_t(B^{s-1}_{p,r})}.
\end{array}
\end{align}
We will temporarily assume that, for any $k\leq n$
\begin{align}
\begin{array}{l}
\|u_k-u_{k-1}\|_{\tilde{L}^{\infty}_t(B^{s-1}_{p,r})}+\|u_k-u_{k-1}\|_{\tilde{L}^2_t(B^s_{p,r})}+
\|h_k-h_{k-1}\|_{\tilde{L}^{\infty}_t(B^{s-1}_{p,r})}
\leq2\times2^{-k}(E_1+E_2).
\end{array}
\end{align}
Then
\begin{align}
\begin{array}{l}
\|u_{n+1}-u_n\|_{\tilde{L}^{\infty}_t(B^{s-1}_{p,r})}+\|u_{n+1}-u_n\|_{\tilde{L}^2_t(B^s_{p,r})}+
\|h_{n+1}-h_n\|_{\tilde{L}^{\infty}_t(B^{s-1}_{p,r})}\\[1ex]
\leq\frac{1}{2}2^{-n}(E_1+E_2)\\[1ex]
+\frac{1}{6}\|u_n-u_{n-1}\|_{\tilde{L}^{\infty}_t(B^{s-1}_{p,r})}+\frac{1}{4}\|u_n-u_{n-1}\|_{\tilde{L}^2_t(B^s_{p,r})}+
\frac{1}{4}\|h_n-h_{n-1}\|_{\tilde{L}^{\infty}_t(B^{s-1}_{p,r})}\\[1ex]
\leq2\times2^{-n-1}(E_1+E_2).
\end{array}
\end{align}
In order to complete the proof of Proposition \ref{p2}, we only need justify the inequalities (3.48) hold for $k=1$, it is obvious that
\begin{align}
\begin{array}{l}
\|u_1-u_0\|_{\tilde{L}^{\infty}_t(B^{s-1}_{p,r})}+\|u_1-u_0\|_{\tilde{L}^2_t(B^s_{p,r})}+
\|h_1-h_0\|_{\tilde{L}^{\infty}_t(B^{s-1}_{p,r})}\\[1ex]
\leq
4(\|u_0\|_{B^s_{p,r}}+\|h_0\|_{B^s_{p,r}})\\[1ex]
\leq E_1+E_2.
\end{array}
\end{align}

This complete the proof of Proposition \ref{p2}.
\subsection{Existence and uniqueness of local solution}
~~~~~In this subsection, we investigate the uniqueness of the local solution to the system (3.1). By Proposition \ref{p2}, the approximative sequence $(u_n,h_n)$ of the problem (2.1) is a Cauchy sequence in $\chi([0,T],s-1,p,r,E_1,E_2)$ with $s>max\{\frac{2}{p},1\}$. So the limit (u,h) is a solution of the Cauchy problem (2.1). From Proposition \ref{p1}, we obtain that this sequence is bounded in $\chi([0,T],s,p,r,E_1,E_2)$. So it's also the Cauchy sequence in $\chi([0,T],s',p,r,E_1,E_2)$ for all $s'<s$ by interpolation. And by Lemma \ref{l7}, the limit is in $\chi([0,T],s,p,r,E_1,E_2)$. Thus we have proved local existence result in Theorem \ref{t1}.

For the uniqueness result in Theorem \ref{t1}, let $(u,h)~~and~~(v,g)$ satisfy the problem (2.1) with the initial data $(u_0,h_0),~(v_0,h_0)\in B^s_{p,r}\times B^s_{p,r}$ respectively. Then we have
\begin{align}
\left\{
\begin{array}{l}
\partial_t(u-v)+u\cdot\nabla(u-v)-\nu\Delta(u-v)=G_1(u,h)-G_1(v,g),\\
\partial_t(h-g)+u\cdot\nabla(h-g)=(u-v)\nabla g+G_2(u,h)-G_2(v,g),\\
(u-v)|_{t=0}=0, (h-g)|_{t=0}=0.
\end{array}
\right.
\end{align}
Using Lemmas \ref{l13}-\ref{l14}, we can get
\begin{align}
\begin{array}{l}
\|u-v\|_{\tilde{L}^{\infty}(B^s_{p,r})}+\|h-g\|_{\tilde{L}^{\infty}(B^s_{p,r})}+\|u-v\|_{\tilde{L}^2(B^{s+1}_{p,r})}\\
\leq C(\|u_0-v_0\|_{B^s_{p,r}}+\|h_0-g_0\|_{B^s_{p,r}})+Ct^{\alpha_1}\|u-v\|_{\tilde{L}^{\infty}(B^s_{p,r})}+
Ct^{\alpha_2}\|h-g\|_{\tilde{L}^{\infty}(B^s_{p,r})}+Ct^{\alpha_3}\|u-v\|_{\tilde{L}^2(B^{s+1}_{p,r})},
\end{array}
\end{align}
here $0<\alpha_1,\alpha_2,\alpha_3<1$.

Then a standard continuous argument gives the uniquness.
\subsection{Continuity}

In this subsection, we will prove that $u,h\in\,\mathcal{C}([0,T];B^s_{p,r})$. First of all, we introduce two useful Lemmas.
\begin{lemm} \label{l33} \cite{H.J}
The time-space estimate for heat equation:

Let $\mathcal{C}$ be an annulus and $\lambda$ a positive real number. Let $u_0,\,f$ satisfy $Supp\,\hat{u}_0,Supp\,\hat{f}(t)\subset \lambda\mathcal{C}$ for all $t\in[0,T]$. Consider u, a solution of
$$\partial_tu-\nu\Delta u=f~~~ and~~~u_{|t=0} =u_0.$$
Then there exists a positive constant $C$, depending only on $\mathcal{C}$, such that for any
$1\leq a\leq b\leq\infty$ and $1\leq p\leq q\leq\infty$, we have
$$\|u\|_{L^q_T(L^b)}\leq C(\nu\lambda^2)^{-\frac{1}{q}}\lambda^{d(\frac{1}{a}-\frac{1}{b})}\|u_0\|_{L^a}+
C(\nu\lambda^2)^{-1+(\frac{1}{p}-\frac{1}{q})}\lambda^{d(\frac{1}{a}-\frac{1}{b})}\|f\|_{L^p_T(L^a)}.$$
\end{lemm}
\begin{lemm}\label{l34} \cite{H.J}

Let $p,p_1,r$ and $s$ be as in the statement of Lemma \ref{l13} with strict inequality in (2.10). Let $f_0\in B^s_{p,r}$, $g\in \tilde{L}^1_T(B^s_{p,r})$, and $v$ be a time-dependent filed such that $v\in L^\rho_T(B^{-M}_{\infty,\infty})$ for some $\rho>1$ and $M>0$, and
$$\nabla v\in L^1_T(B^{\frac{2}{p_1}}_{p_1,\infty}\cap L^\infty),~if ~s<1+\frac{2}{p_1},$$
$$\nabla v\in L^1_T(B^{s-1}_{p_1,\infty}),~if ~s>1+\frac{2}{p_1},~or~s=1+\frac{2}{p_1}~and~r=1.$$
Then the transport equation (2.9) has a unique solution f in

 the space $\mathcal{C}([0,T];B^s_{p,r}),~if~r<\infty$.

 the space $\big(\bigcap_{s'<s}\mathcal{C}([0,T];B^{s'}_{p,r})\big)\bigcap \mathcal{C}_w([0,T];B^s_{p,r}),$ if
  $r=\infty$.

\end{lemm}
\begin{rema}
Although in \cite{H.J}, the continuous conditions request that $g\in L^1_T(B^s_{p,r})$, but by (3.14) in page 134 in \cite{H.J}, $g\in \tilde{L}^1_T(B^s_{p,r})$ does work too.
\end{rema}
From the equations we can get
$u,\,h\in C([0,T];B^{s-2}_{p,r})$, it follows that
$\Delta_ju,\Delta_jh\in \mathcal{C}([0,T];B^s_{p,r})$ for any $j\geq-1$, from which it follows that
$S_ju,S_jh\in \mathcal{C}([0,T];B^s_{p,r})$ for all $j\in N$.
We claim that the sequence of continuous $B^s_{p,r}$-valued fuctions $\{S_ju\}_{j\in N},\{S_jh\}_{j\in N}$ converges uniformly on $[0,T]$. Indeed, by Proposition \ref{l2}, we have
$$\Delta_{j'}(u-S_ju)=\sum_{|j'-j''|\leq1,j''\geq j}\Delta_{j'}\Delta_{j''}u,~~
\Delta_{j'}(h-S_jh)=\sum_{|j'-j''|\leq1,j''\geq j}\Delta_{j'}\Delta_{j''}h,$$
from which it follows that
\begin{align}
\begin{array}{l}
\|u-S_ju\|_{\tilde{L}_T^{\infty}(B^s_{p,r})}\leq C\big(\sum_{j'\geq j-1}2^{j'sr}\|\Delta_{j'}u\|^r_{L_T^{\infty}(L^p)}\big)^{\frac{1}{r}},\\[1ex]
\|h-S_jh\|_{\tilde{L}_T^{\infty}(B^s_{p,r})}\leq C\big(\sum_{j'\geq j-1}2^{j'sr}\|\Delta_{j'}h\|^r_{L_T^{\infty}(L^p)}\big)^{\frac{1}{r}}.
\end{array}
\end{align}
Applying the operator $\Delta_{j'}$ in the first equation of (2.1), we get
\begin{align}
\left\{
\begin{array}{l}
\partial_t\Delta_{j'}u+\Delta_{j'}((u\cdot\nabla)u)-\nu\Delta_{j'}(\Delta u)-\nu\Delta_{j'}(\nu\nabla
(\ln(1+h)\nabla u))+\Delta_{j'}(\nabla h)=0, \\[1ex]
\Delta_{j'}u|_{t=0}=\Delta_{j'}u_0.
\end{array}
\right.
\end{align}
When $j\geq1$, we have the Fourier transform of $\Delta_{j'}u_0$ and $\Delta_{j'}f$  is supported
in an annulus $2^{j'}\mathcal{C}$, by Lemma \ref{l33},  we have
\begin{align}
\begin{array}{l}
\|\Delta_{j'}u\|_{L_T^{\infty}(L^p)}\leq C\|\Delta_{j'}u_0\|_{L^p}+2^{2j'(-1+\frac{1}{2})}\|\Delta_{j'}f\|_{L^2_T(L^p)},
\end{array}
\end{align}
where $f=-\nabla h-u\cdot\nabla u+\nu\nabla(\ln(1+h))\nabla u$.
It follows that
\begin{align}
\begin{array}{l}
\|u-S_ju\|^r_{\tilde{L}_T^{\infty}(B^s_{p,r})}\leq C\sum_{j'\geq j-1}2^{j'sr}\|\Delta_{j'}u\|^r_{L^{\infty}_T(L^p)}\\[1ex]
\leq \sum_{j'\geq j-1}2^{j'sr}\big(\|\Delta_{j'}u_0\|_{L^p}+2^{2j'(-1+\frac{1}{2})}
\|\Delta_{j'}f\|_{L^2_T(L^p)}\big)^r\\[1ex]
\leq C\sum_{j'\geq j-1}2^{j'sr}\|\Delta_{j'}u_0\|^r_{L^p}+C\sum_{j'\geq j-1}\big(
2^{j'(s-1)}\|\Delta_{j'}f\|_{L^2_T(L^p)}\big)^r
\end{array}
\end{align}

The first term clearly tends to $0$ when $j$ goes to $\infty$. For the second term, if $f\in \tilde{L}^2_T(B^{s-1}_{p,r})$, we also have that it tends to $0$ when $j$ goes to $\infty$.
Clearly, by $h\in \tilde{L}^\infty_T(B^s_{p,r})$, we have that
\begin{align}
\begin{array}{l}
\|\nabla h\|_{\tilde{L}^2_T(B^{s-1}_{p,r})}\leq CT^{\frac{1}{2}}\|h\|_{\tilde{L}^\infty_T(B^s_{p,r})}.
\end{array}
\end{align}
For $u\cdot\nabla u$, by Lemmas \ref{l10}, \ref{l11}. we have
\begin{align}
\begin{array}{l}
\|u\cdot\nabla u\|_{\tilde{L}^2_T(B^{s-1}_{p,r})}\\[1ex]
\leq C\|u\|_{L^\infty_T(L^\infty)}\|\nabla u\|_{\tilde{L}^2_T(B^{s-1}_{p,r})}+C\|\nabla u\|_{\tilde{L}^2_T(L^\infty)}\|u\|_{\tilde{L}^\infty_T(B^{s-1}_{p,r})}\\[1ex]
+C\|u\|_{\tilde{L}^\infty_T(B^s_{p,r})}\|\nabla u\|_{L^2_T(B^{-1}_{\infty,\infty})}\\[1ex]
\leq C\|u\|_{\tilde{L}^\infty_T(B^s_{p,r})}\|u\|_{\tilde{L}^2_T(B^{s+1}_{p,r})}.
\end{array}
\end{align}
As regards $\nabla(\ln(1+h))\nabla u$, still by Lemmas \ref{l10},\ref{l11}, we can get
\begin{align}
\begin{array}{l}
\|\nabla(\ln(1+h))\nabla u\|_{\tilde{L}^2_T(B^{s-1}_{p,r})}\\[1ex]
\leq C\|\nabla(\ln(1+h))\|_{\tilde{L}^\infty_T(B^{s-1-\frac{2}{p}}_{\infty,\infty})}\|\nabla u\|_{L^2_T(B^{\frac{2}{p}}_{p,r})}+C\|\nabla u\|_{L^2_T(L^\infty)}\|\nabla(\ln(1+h))\|_{\tilde{L}^\infty_T(B^{s-1}_{p,r})}\\[1ex]
+C\|\nabla(\ln(1+h))\|_{\tilde{L}^\infty_T(B^{s-1}_{p,r})}\|\nabla u\|_{L^2_T(B^0_{\infty,\infty})}\\[1ex]
\leq C\|u\|_{\tilde{L}^\infty_T(B^s_{p,r})}\|u\|_{\tilde{L}^2_T(B^{s+1}_{p,r})}.
\end{array}
\end{align}
 This completes the proof of continuity for $u$.

As regards $h$, by the second equation of (2.1) and Lemma \ref{l34}, it suffices to prove that
$div\,u,~hdiv\,u\in \tilde{L}^1_T(B^s_{p,r})$.
It's obvious that
\begin{align}
\begin{array}{l}
\|div\,u\|_{\tilde{L}^1_T(B^s_{p,r})}\leq CT^{\frac{1}{2}}\|u\|_{\tilde{L}^2_T(B^s_{p,r})}.
\end{array}
\end{align}
By Lemma \ref{l12}, we have
\begin{align}
\begin{array}{l}
\|h\,div\,u\|_{\tilde{L}^1_T(B^s_{p,r})}\\[1ex]
\leq C\|h\|_{L^\infty_T(L^\infty)}\|div\,u\|_{\tilde{L}^1_T(B^s_{p,r})}+C\|div\,u\|_{L^1_T(L^\infty)}
\|h\|_{\tilde{L}^\infty_T(B^s_{p,r})}\\[1ex]
\leq C\|u\|_{\tilde{L}^\infty_T(B^s_{p,r})}\|u\|_{\tilde{L}^2_T(B^{s+1}_{p,r})}.
\end{array}
\end{align}
This completes the proof of continuty for $h$.

\section{Global existence}
In this section, we will get the global existence of the
system with small enough initial data $(u_0, h_0)$ in
$B^s_{p,r}\times B^s_{p,r}$, $1\leq p\leq2$, $1\leq r<\infty$ and $\frac{2}{p}<s<\frac{2}{p}+1$(the case of $s\geq 1+\frac{2}{p}$ is easier).

By the imbedding theorem, it's obvious that $(u_0,h_0)$ in $H^{s-\frac{2}{p}+1-\ep}\times H^{s-\frac{2}{p}+1-\ep}$, and
$$\|u_0\|_{H^{s-\frac{2}{p}+1-\ep}}\leq C\|u_0\|_{B_{p,r}^s},~\|h_0\|_{H^{s-\frac{2}{p}+1-\ep}}\leq C\|h_0\|_{B_{p,r}^s},$$
where $\ep$ is small enough such that $\frac{2}{p}<s-\ep<\frac{2}{p}+1$ too. For the sake of convenience, letting $s_1=s-\frac{2}{p}+1-\ep$, thus $1<s_1<2$.

Then by Lemma \ref{h} there exits a positive real number $\eta>0$ such that, if $$\|u_0\|_{B_{p,r}^s}+\|h_0\|_{B_{p,r}^s}\leq \eta,$$
the system (2.1) have a unique global solution $(u,h)$ in $\tilde{L}^\infty_T(H^{s_1})\cap
\tilde{L}^2_T(H^{s_1+1})\times \tilde{L}^\infty_T(H^{s_1})$.  Thus it's sufficient to prove the following Proposition.
\begin{prop}\label{p4}
Let $(u_0,~h_0)\in B_{p,r}^s\times B_{p,r}^s$, $s>\frac{2}{p},~1\leq p\leq2$,  $1\leq r<\infty$,
$\|u_0\|_{B_{p,r}^s}+\|h_0\|_{B_{p,r}^s}\leq \eta$ such that system(2.1) has
 a unique global solution $(u,h)$ satisfy for any $T\geq0$,
$$u\in \tilde{L}^{\infty}([0,T];H^{s_1})\cap \tilde{L}^2([0,T];H^{s_1+1})$$
$$h\in \tilde{L}^{\infty}([0,T];H^{s_1}).$$
Then we claim that for any $T\geq0$, $\|u(T)\|_{B_{p,r}^s},\|h(T)\|_{B_{p,r}^s}$ is finite.
\end{prop}
Proof: First of all, by $u\in \tilde{L}^{\infty}([0,T];H^{s_1})\cap \tilde{L}^2([0,T];H^{s_1+1})$ and the interpolation inequality, we have
$$u\in \tilde{L}^{\rho_1}([0,T];H^{\frac{s_1+3}{2}}), ~where~ \frac{\theta}{s_1}+
 \frac{1-\theta}{s_1+1}=\frac{2}{s_1+3},~\frac{\theta}{\infty}+\frac{1-\theta}{2}=\frac{1}{\rho_1} $$
 $\theta\in(0,1)$, $\rho_1\in(2,\infty)$.

For $\|h(T)\|_{B^s_{p,r}}$, by Lemma \ref{l13}, we have
\begin{align}
\begin{array}{l}
\|h\|_{\tilde{L}^\infty_T(B^s_{p,r})}\leq exp\big(C_0\int_0^T\|\nabla u\|_{B^1_{2,\infty}\cap L^\infty}dt\big)\big(\|h_0\|_{B^s_{p,r}}+\|div\,u\|_{\tilde{L}^1_T(B^s_{p,r})}+\|h\,div\,u\|_{\tilde{L}^1_T(B^s_{p,r})}\big).
\end{array}
\end{align}
It's obvious that
\begin{align}
\begin{array}{l}
exp\big(C_0\int_0^T\|\nabla u\|_{B^1_{2,\infty}\cap L^\infty}dt\big)\\[1ex]
\leq exp\big(2C_0T^{\frac{1}{2}}\|u\|_{L^2_T(H^{s_1+1})}\big)\\[1ex]\
\leq C
\end{array}
\end{align}
and
\begin{align}
\begin{array}{l}
\|div\,u\|_{\tilde{L}^1_T(B^s_{p,r})}\leq\|u\|_{L^1_T(B^{s+1}_{p,r})}.
\end{array}
\end{align}
By Lemmas \ref{l10}, \ref{l11}, we have
\begin{align}
\begin{array}{l}
\|h\,div\,u\|_{\tilde{L}^1_T(B^s_{p,r})}\\[1ex]
\leq \|T_hdiv\,u\|_{\tilde{L}^1_T(B^s_{p,r})}+\|T_{div\,u} h\|_{\tilde{L}^1_T(B^s_{p,r})}+
\|R(h,div\,u)\|_{\tilde{L}^1_T(B^s_{p,r})}   \\[1ex]
\leq C\|h\|_{L^\infty_T(L^\infty)}\|div\,u\|_{\tilde{L}^1_T(B^s_{p,r})}+C\|T_{div\,u} h\|_{L^1_T(B^s_{p,r})}+C\|h\|_{L^\infty_T(B^0_{\infty,\infty})}\|div\,u\|_{\tilde{L}^1_T(B^s_{p,r})}\\[1ex]
\leq C\|h\|_{L^\infty_T(H^{s_1})}\|u\|_{\tilde{L}^1_T(B^{s+1}_{p,r})}+C\|div\, u\|_{L^2_T(L^\infty)}\|h\|_{L^2_T(B^s_{p,r})}\\[1ex]
\leq C\|h\|_{L^\infty_T(H^{s_1})}\|u\|_{\tilde{L}^1_T(B^{s+1}_{p,r})}+C\| u\|_{L^2_T(H^{s_1+1})}\|h\|_{L^2_T(B^s_{p,r})}.
\end{array}
\end{align}
Combining (4.1)-(4.4), we get
\begin{align}
\begin{array}{l}
\|h\|_{\tilde{L}^\infty_T(B^s_{p,r})}\leq C\|h_0\|_{B^s_{p,r}}+C(1+\|h\|_{L^\infty_T(H^{s_1})})\|u\|_{\tilde{L}^1_T(B^{s+1}_{p,r})}+C\| u\|_{L^2_T(H^{s_1+1})}\|h\|_{L^2_T(B^s_{p,r})}\\[1ex]
\leq C+C\|u\|_{\tilde{L}^2_T(B^{s+1}_{p,r})}+C\|h\|_{L^2_T(B^s_{p,r})}.
\end{array}
\end{align}
Here and hereafter, $C$ is a constant depending on $T,s,p,r,\|u_0\|_{B^s_{p,r}},\|h_0\|_{B^s_{p,r}},\nu$.

Then by Lemmas \ref{l14},  we get
\begin{align}
\begin{array}{l}
\|u\|_{\tilde{L}^\infty_T(B^s_{p,r})}\leq  C_0exp\big(\int_0^T\|\nabla u\|_{B^1_{2,\infty}\cap L^\infty}dt\big)\times\\[1ex]
 \big(\|u_0\|_{B^s_{p,r}}+C\|\nabla h\|_{\tilde{L}^\rho_T(B^{s-2+\frac{2}{\rho}}_{p,r})}+C\|\nabla(\ln(1+h))\nabla u\|_{\tilde{L}^\rho_T(B^{s-2+\frac{2}{\rho}}_{p,r})}\big)\\[1ex]
\leq C+C\|h\|_{L^\rho_T(B^s_{p,r})}+C\|\nabla(\ln(1+h))\nabla u\|_{\tilde{L}^\rho_T(B^{s-2+\frac{2}{\rho}}_{p,r})},
\end{array}
\end{align}
here $\rho\in(2,\rho_1)$.
Then by Lemma \ref{l10} and corollary \ref{r}, we have
\begin{align}
\begin{array}{l}
\|\nabla(\ln(1+h))\nabla u\|_{\tilde{L}^\rho_T(B^{s-2+\frac{2}{\rho}}_{p,r})}\\[1ex]
\leq C\|\nabla u\|_{L^{\rho_1}_T(L^\infty)}\|\nabla(\ln(1+h))\|_{\tilde{L}^{\rho_2}_T(B^{s-2+\frac{2}{\rho}}_{p,r})}+
C\|\nabla(\ln(1+h))\|_{L^\infty_T(B^{s_1-2}_{\infty,\infty})}\|\nabla u\|_{\tilde{L}^\rho_T(B^{s+\frac{2}{\rho}-s_1}_{p,r})}\\[1ex]
\leq C\|u\|_{L^{\rho_1}_T(H^{\frac{s_1+3}{2}})}\|h\|_{L^{\rho_2}_T(B^s_{p,r})}+
C\|h\|_{L^\infty_T(H^{s_1})}\|u\|_{\tilde{L}^\rho_T(B^{s+\frac{2}{\rho}-s_1+1}_{p,r})}.
\end{array}
\end{align}
here $\rho,\rho_1,\rho_2$ satisfy
 $$\frac{1}{\rho_1}+\frac{1}{\rho_2}=\frac{1}{\rho}.$$

Combining (4.6)-(4.7), we can obtain
\begin{align}
\begin{array}{l}
\|u\|_{\tilde{L}^\infty_T(B^s_{p,r})}\\[1ex]
\leq C+C\|h\|_{L^\rho_T(B^s_{p,r})}+C\|h\|_{L^{\rho_2}_T(B^s_{p,r})}+
C\|u\|_{\tilde{L}^\rho_T(B^{s+\frac{2}{\rho}-s_1+1}_{p,r})}\\[1ex]
\leq C+C\|h\|_{L^{\rho_2}_T(B^s_{p,r})}+C\|u\|_{\tilde{L}^\rho_T(B^{s+\frac{2}{\rho}-s_1+1}_{p,r})}.
\end{array}
\end{align}

Similarly with $\|u\|_{\tilde{L}^\infty_T(B^s_{p,r})}$, using Lemma \ref{l10} and Corollary \ref{r} again, we have
\begin{align}
\begin{array}{l}
\|u\|_{\tilde{L}^\rho_T(B^{s+\frac{2}{\rho}-s_1+1}_{p,r})}\leq C\big(\|u_0\|_{B^{s-s_1+1}_{p,r}}+C\|\nabla h+\nabla(\ln(1+h))\nabla u\|_{\tilde{L}^\rho_T(B^{s+\frac{2}{\rho}-s_1-1}_{p,r})}\big)\\[1ex]
\leq C+C\|h\|_{L^\rho_T(B^s_{p,r})}+C\|\nabla(\ln(1+h))\nabla u\|_{\tilde{L}^\rho_T(B^{s+\frac{2}{\rho}-s_1-1}_{p,r})}\\[1ex]
\leq C+C\|h\|_{L^\rho_T(B^s_{p,r})}+C\|\nabla u\|_{L^{\rho_1}_T(L^\infty)}\|\nabla(\ln(1+h))\|_{\tilde{L}_T^{\rho_2}(B^{s+\frac{2}{\rho}-s_1-1}_{p,r})}\\[1ex]
+C\|\nabla(\ln(1+h))\|_{L^\infty_T(B^{s_1-2}_{\infty,\infty})}\|\nabla u\|_{\tilde{L}^\rho_T(B^{s+\frac{2}{\rho}-2s_1+1}_{p,r})}\\[1ex]
\leq C+C\|h\|_{L^{\rho_2}_T(B^s_{p,r})}+C\|u\|_{\tilde{L}^\rho_T(B^{s+\frac{2}{\rho}-2s_1+2}_{p,r})},
\end{array}
\end{align}

here we use the fact that $$s+\frac{2}{\rho}-s_1-1=s+\frac{2}{\rho}-(s-\frac{2}{p}+1-\ep)-1=\frac{2}{\rho}+\ep+\frac{2}{p}-2>\frac{2}{p}-2.$$

Combining (4.8)-(4.9), we can obtain
\begin{align}
\begin{array}{l}
\|u\|_{\tilde{L}^\infty_T(B^s_{p,r})}\\[1ex]
\leq C+C\|h\|_{L^{\rho_2}_T(B^s_{p,r})}+C\|u\|_{\tilde{L}^\rho_T(B^{s+\frac{2}{\rho}-2s_1+2}_{p,r})}.
\end{array}
\end{align}
Then by Lemma \ref{l10}, we get

\begin{align}
\begin{array}{l}
\|u\|_{\tilde{L}^\rho_T(B^{s+\frac{2}{\rho}-2s_1+2}_{p,r})}\leq C\big(\|u_0\|_{B^{s-2s_1+2}_{p,r}}+C\|\nabla h+\nabla(\ln(1+h))\nabla u\|_{\tilde{L}^\rho_T(B^{s+\frac{2}{\rho}-2s_1}_{p,r})}\big)\\[1ex]
\leq C+C\|h\|_{L^\rho_T(B^s_{p,r})}+C\|\nabla(\ln(1+h))\nabla u\|_{\tilde{L}^\rho_T(B^{s+\frac{2}{\rho}-2s_1}_{p,r})}\\[1ex]
\leq C+C\|h\|_{L^\rho_T(B^s_{p,r})}+C\|T_{\nabla u}\nabla(\ln(1+h))\|_{\tilde{L}^\rho_T(B^{s+\frac{2}{\rho}-2s_1}_{p,r})}\\[1ex]+C\|T_{\nabla(\ln(1+h))}\nabla u\|_{\tilde{L}^\rho_T(B^{s+\frac{2}{\rho}-2s_1}_{p,r})}+\|R(\nabla(\ln(1+h)),\nabla u)\|_{\tilde{L}^\rho_T(B^{s+\frac{2}{\rho}-2s_1}_{p,r})}\\[1ex]
\leq C+C\|h\|_{L^\rho_T(B^s_{p,r})}+C\|\nabla u\|_{L^{\rho_1}_T(L^\infty)}\|\nabla(\ln(1+h))\|_{\tilde{L}_T^{\rho_2}(B^{s+\frac{2}{\rho}-2s_1}_{p,r})}\\[1ex]
+C\|\nabla(\ln(1+h))\|_{L^\infty_T(B^{s_1-2}_{\infty,\infty})}\|\nabla u\|_{\tilde{L}^\rho_T(B^{s+\frac{2}{\rho}-3s_1+2}_{p,r})}+\|R(\nabla(\ln(1+h)),\nabla u)\|_{\tilde{L}^\rho_T(B^{s+\frac{2}{\rho}-2s_1}_{p,r})}\\[1ex]
\leq C+C\|h\|_{L^{\rho_2}_T(B^s_{p,r})}+C\|u\|_{\tilde{L}^\rho_T(B^{s+\frac{2}{\rho}-3s_1+3}_{p,r})}+\|R(\nabla(\ln(1+h)),\nabla u)\|_{\tilde{L}^\rho_T(B^{s+\frac{2}{\rho}-2s_1}_{p,r})}.
\end{array}
\end{align}

For The remainder term $\|R\big(\nabla u,\nabla(\ln(1+h))\big)\|_{\tilde{L}^\rho_T(B_{1,r}^{s+\frac{2}{\rho}-2s_1})}$, if
$s+\frac{2}{\rho}-2s_1>\frac{2}{p}-2$, by Corollary \ref{r} and (4.11), we have
\begin{align}
\begin{array}{l}
\|R\big(\nabla u,\nabla(\ln(1+h))\big)\|_{\tilde{L}^\rho_T(B_{p,r}^{s+\frac{2}{\rho}-2s_1})}\leq C
\|R\big(\nabla u,\nabla(\ln(1+h))\big)\|_{\tilde{L}^\rho_T(B_{1,r}^{s+\frac{2}{\rho}-2s_1+2-\frac{2}{p}})}\\[1ex]
\leq C\|\nabla (\ln(1+h))\|_{L^\infty_T(B^{s_1-\frac{2}{p}}_{p',\infty})}
\|\nabla u\|_{\tilde{L}^\rho_T(B^{s+\frac{2}{\rho}-3s_1+2}_{p,r})}\\[1ex]
\leq C\|h\|_{L_T^\infty(H^{s_1})}\|u\|_{\tilde{L}^\rho_T(B^{s+\frac{2}{\rho}-3s_1+3}_{p,r})}.
\end{array}
\end{align}

If $s+\frac{2}{\rho}-2s_1\leq\frac{2}{p}-2$, we have

\begin{align}
\begin{array}{l}
\|R\big(\nabla u,\nabla(\ln(1+h))\big)\|_{\tilde{L}^\rho_T(B_{p,r}^{s+\frac{2}{\rho}-2s_1})}\\[1ex]
\leq C
\|R\big(\nabla u,\nabla(\ln(1+h))\big)\|_{\tilde{L}^\rho_T(B_{1,r}^\ep)}\\[1ex]\
\leq C\|\nabla(\ln(1+h))\|_{L^\infty_T(B^{s_1-1}_{2,\infty})}\|\nabla u\|_{\tilde{L}^\rho_T(B^{\ep+1-s_1}_{2,r})}\\[1ex]
\leq C\|h\|_{\tilde{L}^\infty_T(H^{s_1})}\|u\|_{L^\rho_T(B_{2,r}^{\ep+2-s_1})}\\[1ex]
\leq C\|h\|_{\tilde{L}^\infty_T(H^{s_1})}\|u\|_{L^{\rho_1}_T(H^{\frac{s_1+3}{2}})}\\[1ex]
\leq C.
\end{array}
\end{align}
where $0<\ep<<1$.

Combining (4.10)-(4.13), we have
\begin{align}
\begin{array}{l}
\|u\|_{\tilde{L}^\infty_T(B^s_{p,r})}\\[1ex]
\leq C+C\|h\|_{L^{\rho_2}_T(B^s_{p,r})}+C\|u\|_{\tilde{L}^\rho_T(B^{s+\frac{2}{\rho}-3s_1+3}_{p,r})}.
\end{array}
\end{align}

Because of $s_1>1$, there exists a positive integer $n$(bigger than 3), such that $s_1>1+\frac{1}{n}\frac{2}{\rho}$.

Repeating the calculation $n-3$ times ,we have
\begin{align}
\begin{array}{l}
\|u\|_{\tilde{L}^\infty_T(B^s_{p,r})}\\[1ex]
\leq C+C\|h\|_{L^{\rho_2}_T(B^s_{p,r})}+C\|u\|_{\tilde{L}^\rho_T(B^{s+\frac{2}{\rho}-ns_1+n}_{p,r})}\\[1ex]
\leq C+C\|h\|_{L^{\rho_2}_T(B^s_{p,r})}+C\|u\|_{L^\rho_T(B^s_{p,r})}.
\end{array}
\end{align}
Similarly for $\|u\|_{\tilde{L}^2_T(B^{s+1}_{p,r})}$, we have the same estimate

\begin{align}
\begin{array}{l}
\|u\|_{\tilde{L}^2_T(B^{s+1}_{p,r})}\\[1ex]
\leq C+C\|h\|_{L^{\rho_2}_T(B^s_{p,r})}+C\|u\|_{L^\rho_T(B^s_{p,r})}.
\end{array}
\end{align}

thus combining (4.5),(4.15),(4.16), we have
\begin{align}
\begin{array}{l}
\|h\|_{\tilde{L}^\infty_T(B^s_{p,r})}+\|u\|_{\tilde{L}^\infty_T(B^s_{p,r})}\\[1ex]
\leq C+C\|h\|_{L^{\rho_2}_T(B^s_{p,r})}+C\|u\|_{L^\rho_T(B^s_{p,r})}\\[1ex]
\leq C+C\|h\|_{L^{\rho_2}_T(B^s_{p,r})}+C\|u\|_{L^{\rho_2}_T(B^s_{p,r})}.
\end{array}
\end{align}
Taking $\rho_2$ power in both sides, we get
\begin{align}
\begin{array}{l}
\big(\|h\|_{\tilde{L}^\infty_T(B^s_{p,r})}+\|u\|_{\tilde{L}^\infty_T(B^s_{p,r})}\big)^{\rho_2}\\[1ex]
\leq C+C\int^T_0\big(\|h(t)\|_{(B^s_{p,r})}+C\|u(t)\|_{(B^s_{p,r})}\big)^{\rho_2}dt.
\end{array}
\end{align}
In the view of Gronwall inequality, we have
\begin{align}
\begin{array}{l}
\big(\|h\|_{\tilde{L}^\infty_T(B^s_{p,r})}+\|u\|_{\tilde{L}^\infty_T(B^s_{p,r})}\big)^{\rho_2}\\[1ex]
\leq Ce^{CT}.
\end{array}
\end{align}

This completes the prove of Proposition \ref{p4}.

\vspace*{2em} \noindent\textbf{Acknowledgements} This work was
partially supported by NNSFC (No.11271382 No. 10971235), RFDP (No.
20120171110014), and the key project of Sun Yat-sen University.

\phantomsection
\addcontentsline{toc}{section}{\refname}
%Ìí¼Ó²Î¿¼ÎÄÏ×µ½ÊéÇ©£¬ºê°ü hyperref

\end{document}